\documentclass[titlepage,draft,12pt]{article} 
\usepackage{amsfonts,latexsym,pstricks,pst-plot,graphicx} 
\usepackage[a4paper]{geometry}

\title{Asymptotic limits for mildly degenerate
Kirchhoff equations} 
\author{Marina Ghisi\vs\\ {\normalsize
Universit\`a degli Studi di Pisa} \\{\normalsize Dipartimento di
Matematica ``Leonida Tonelli''}\\
{\normalsize 
PISA (Italy)}\\  
{\normalsize e-mail: \texttt{ghisi@dm.unipi.it}}
}
\newcommand{\ep}{\varepsilon}
\newcommand{\re}{{\mathbb{R}}}

\newcommand{\vs}{\vspace{.2cm}}

\newcommand{\qed}{{\penalty 10000\mbox{$\quad\Box$}}\bigskip}

\newcommand{\foralll}{\forall\:}

\newtheorem{thmbibl}{Theorem}

\newcommand{\m}[1]{\au{#1}^{2\gamma}}
\newcommand{\au}[1]{|A^{1/2}#1|}
\newcommand{\auq}[1]{|A^{1/2}#1|^{2}}

\newcommand{\da}{D(A)}
\newcommand{\dau}{D(A^{1/2})}

\newcommand{\eB}[1]{e^{2#1 B_{\ep}(t)}}

\newcommand{\uep}{u_{\ep}}

\newtheorem{thm}{Theorem}[section]
\newtheorem{rmk}[thm]{Remark}
\newtheorem{prop}[thm]{Proposition}

\newtheorem{lemma}[thm]{Lemma}

\begin{document}
\maketitle
\begin{abstract}
	We consider the second order Cauchy problem
	$$\ep\uep''+\m{\uep}A\uep+\uep'=0, \hspace{3em} 
	\uep(0)=u_0\neq 0,
	\hspace{2em}\uep'(0)=u_1$$
	where $\ep>0$, $H$ is an Hilbert space, $A$ is a self-adjoint
	positive operator on $H$ with dense domain $D(A)$,
	$(u_{0},u_{1})\in\da\times\dau$, and  $\gamma > 0$.
	
	We study accurately the  decay as $t$ goes 
	to infinity of the 
	solutions, provided that $\ep$
	is small enough. In particular we obtain a new estimate on 
	$\uep''$ and we show that the renormalized functions 
	$(1+t)^{1/(2\gamma)}(\uep(t),(1+t)\uep'(t))$ 
	have a non zero limit in $D(A)\times D(A^{1/2})$ as $t$ goes 
	to infinity.  
	Moreover we  calculate 
	explicitly the norm of these limit and we  prove that they do
	not depend on the initial 
	data.

\vspace{1cm}

\noindent{\bf Mathematics Subject Classification 2000 (MSC2000):}
35B25, 35B40, 35L80.

\vspace{1cm} 
\noindent{\bf Key words:}  degenerate
damped hyperbolic equations, Kirchhoff
equations, decay rate of solutions.
\end{abstract}
 
\section{Introduction}

Let $H$ be a real Hilbert space.  Given $x$ and $y$ in $H$, $|x|$
denotes the norm of $x$, and $\langle x,y\rangle$ denotes the scalar
product of $x$ and $y$.  Let $A$ be a self-adjoint linear operator on
$H$ with dense domain $D(A)$.  We always assume that $A$ is
coercive, namely $\langle Au,u\rangle\geq \sigma_{0}|u|^{2}$ for every $u\in D(A)$.
For any such operator the power $A^{\alpha}$ is defined for every
$\alpha\geq 0$ in a suitable domain $D(A^{\alpha})$. 

For every $\ep>0$ we consider the second order Cauchy problem
\begin{eqnarray}
	 &  & \ep\uep''(t)+\m{\uep(t)}A\uep(t)+\uep'(t)=0, \hspace{2em}
	\foralll t\geq 0,
	\label{pbm:h-eq}  \\
	 &  & \uep(0)=u_0\neq 0, \hspace{2em}\uep'(0)=u_1,
	\label{pbm:h-data}
\end{eqnarray}
 where $u_{0}\in D(A)$, $u_{1}\in D(A^{1/2})$. This problem is just an abstract
setting (with $m(r) = r^{\gamma}$) of the initial boundary value
problem for the hyperbolic partial differential equation (PDE)
\begin{equation}
	\ep u^{\ep}_{tt}(t,x)-
	m{\left(\int_{\Omega}\left|\nabla u^\ep(t,x)\right|^2\,dx\right)}
	\Delta u^\ep(t,x)+u^\ep_t(t,x)=0  \\
	\label{pbm:h-concr}
\end{equation}
in a bounded open set $\Omega\subseteq\re^{n}$, where $m$ is a non 
negative function.  This equation is a model for
the damped small transversal vibrations of an elastic string ($n=1$)
or membrane ($n=2$) with uniform density $\ep$.

Equations as (\ref{pbm:h-eq}) or (\ref{pbm:h-concr}) have been 
intensely studied from `80 both in the case of an operator $A$ 
coercive and in the case of an only nonnegative operator. In 
particular it was extensively considered the case of a function $m$ 
in the $C^{1}$ class, both in the 
nondegenerate case ($m\geq c >0$) and in the mildly degenerate case 
($m(|A^{1/2}u_{0}|^{2})\neq 0$). For a more complete discussion 
on this argument we refer to the survey \cite{trieste} and to the references  
contained therein. Here we concentrate shortly only on few  of the results 
concerning the existence of global solutions and their behaviour at 
the infinity. First of all let us remind that there are not 
substantial differences between coercive and only nonnegative 
operators with respect to the existence of global solutions, while 
there are differences regarding the asymptotic behaviour. Indeed in 
the case of only nonnegative operators the estimates that one can 
obtain on $|A^{1/2}\uep|$ are in general worse (see \cite{gg:k-decay}).  All the 
results we state explicitly regarding the decay of the solutions  must then be thought in the coercive 
case.  We point out that the existence of a global solution 
for small data was proved firstly in the nondegenerate case of $m\geq 
c > 0$ in
\cite{debrito1} and in \cite{y}, then in the mildly 
degenerate case of (\ref{pbm:h-eq}) - (\ref{pbm:h-data})  by 
\textsc{K. Nishihara and Y. Yamada} \cite{ny} if $\gamma \geq 1$ and 
in \cite{jmaa} - \cite{g:non-lip} when $0< \gamma 
< 1$. In all these papers it was also considered the behaviour at the 
infinity of the solutions and some, in general non 
optimal, estimates were obtained  (see also \cite{debrito2}, 
\cite{yamazaki}, \cite{hy}
\cite{n-natma} and  \cite{ono-caa}  for 
the nondegenerate case). In the mildly degenerate case of 
(\ref{pbm:h-eq}) good 
estimates on $|A^{1/2}\uep|$, $|A\uep|$, $|\uep'|$ were proved firstly by  
\textsc {T. Mizumachi} (\cite{mizu-ade}, \cite{mizu-nc}) and \textsc{K. Ono}
(\cite{ono-kyushu}, \cite{ono-aa}) when $\gamma =1$ and 
then for any $\gamma > 0$ in \cite{gg:k-decay}, \cite{jde2} (see 
Theorem \ref{A} for the precise statement in the coercive case). It is clear that 
estimates on $|A^{1/2}\uep|$, $|A\uep|$, $|\uep'|$ produce estimates also 
 on $\uep''$. These estimates are in general not sharp, as 
it was shown when $\gamma = 1$ in  \cite{ono-kyushu}. 
Indeed in this last  case the obvious estimate gives 
$(1+t)^{3}|\uep''(t)|^{2} \leq C_{\ep}$, while K. Ono proved  that at 
least one has $(1+t)^{4}|\uep''(t)|^{2} \leq C_{\ep}$. 

We proposed to clarify more precisely the behaviour of the solutions 
of  (\ref{pbm:h-eq}) - (\ref{pbm:h-data}), by studying the problem 
from a new point of view. We prove indeed that not only estimates as 
(\ref{h1})-(\ref{h12}) hold true, but in fact the renormalized 
functions $(1+t)^{1/(2\gamma)}(\uep(t),(1+t)\uep'(t))$ have a \emph{non 
zero} limit 
$(u_{\ep,\infty}, v_{\ep,\infty})$ in $D(A)\times D(A^{1/2})$. Of 
such limits we can calculate explicitly the norms, that do not 
depend on the initial data or $\ep$. As a consequence we obtain sharp 
estimates on the decay of the solutions of (\ref{pbm:h-eq}) 
in $D(A)\times D(A^{1/2})$.  Equally important we can express 
explicitly the relation between  $u_{\ep,\infty}$ and 
$v_{\ep,\infty}$.
This allows us also to 
obtain better estimates on $\uep''$, indeed for example  in the case of 
$\gamma = 1$ our estimates give (see Theorem \ref{t2}):  $$(1+t)^{5}|\uep''(t)|^{2} \leq 
C_{\ep}.$$  
For our purpose in the following we assume that $H$ has a countable basis made by 
eigenvectors of $A$, that is obviously verified  in the concrete case 
of (\ref{pbm:h-concr}). It is well
known that in fact it is enough to assume that  initial data in 
(\ref{pbm:h-data}) can be written in Fourier series with respect to 
 eigenvectors of $A$, because in such a case this is also true 
for the solution of (\ref{pbm:h-eq}). This hypothesis allows us to 
prove (see Theorem \ref{t2} and (\ref{B32b}), (\ref{B4b})) that actually  
the norms of $u_{\ep,\infty}$ in $D(A)$ and of $v_{\ep,\infty}$ in 
$D(A^{1/2})$  depend only on 
$\gamma$ and on the smallest eigenvalue of $A$. 
This surprising behaviour depends on the  fact that actually 
the behaviour of the solutions of  (\ref{pbm:h-eq}) is due only to 
the components of  $\uep$ related to the smallest eigenvalue of $A$  
(see Theorem \ref{t1}). This type of estimates so precise should 
allow us to obtain decay-error estimates in the study of the singular 
perturbation problem,  that consists in setting formally $\ep = 0$ in 
(\ref{pbm:h-eq}) and study the difference between $\uep$ and the 
solution of the new first order problem obtained in such a way 
(see \cite{trieste} and \cite{yamazaki} for an 
introduction of the problem and its treatment in the non degenerate 
case).

The outline of the paper is the following. In Section 
\ref{sec:statements} to begin with we recall the hold result we need 
on the existence of global solutions of (\ref{pbm:h-eq}) and their 
decay and 
we introduce some preliminary notations, then we state the main 
results. In section \ref{sec:proofs} we prove the results. This last 
section is divided in various parts. First of all we prove a general linear 
result that we use in particular for proving Theorem \ref{t1}, then 
we study the properties of the components of the solution of 
(\ref{pbm:h-eq}) and finally prove Theorem \ref{t2}. Let us stress 
that  Theorem \ref{t1} is in fact a linear result whereas  proof of Theorem \ref{t2} requires a new type of nonlinear approach.

\setcounter{equation}{0}

\section{Statements}\label{sec:statements}

\subsection{Notations and preliminaries}
Let us stress that we assumed that the operator $A$ is coercive. The following result is well known and it is a consequence of
\cite{ny}, \cite{g:non-lip} (see also \cite{gg:k-dissipative} for the 
study of the  
case of more general functions $m$) for the part concerning the
existence of global solutions, while it follows from \cite{gg:k-decay} \cite{jde2} (see also
 \cite{mizu-ade}, \cite{mizu-nc}, \cite{ny},
\cite{ono-kyushu}, \cite{ono-aa}) for the part concerning the  decay of solutions.

\begin{thmbibl} \label{A} \em{ Let $(u_{0},u_{1})\in D(A)\times
D(A^{1/2})$.  Let  $\gamma > 0$, then for $\ep$ small the mildly
degenerate problem (\ref{pbm:h-eq}), (\ref{pbm:h-data}) has a unique
global solution $$\uep\in C^{2}([0,+\infty[,H) \cap  
C^{1}([0,+\infty[,D(A^{1/2})) \cap  C^{0}([0,+\infty[,D(A))$$ such that
       \begin{equation}
	   \frac{K_{1}}{(1+t)^{1/\gamma}}\leq \auq{\uep(t)}\leq
	   \frac{K_{2}}{(1+t)^{1/\gamma}}\label{h1} \quad\quad\forall
	   t\geq 0,  
       \end{equation}
		\begin{equation}
		\frac{K_{1}}{(1+t)^{1/\gamma}}\leq |A\uep(t)|^{2}\leq
		\frac{K_{2}}{(1+t)^{1/\gamma}} \quad\quad\forall t\geq
		0,
	   \label{h11}
       \end{equation}
      \begin{equation} |\uep'(t)|^{2}\leq
      \frac{K_{2}}{(1+t)^{2+1/\gamma}} \quad\quad\forall t\geq 0;
		 \label{h12}
		 \end{equation}
		 where the constants $K_{1}$ and $K_{2}$ do not depend
		 on $\ep$.}
		\end{thmbibl}
	In the following we assume always that $\ep\leq 1$ is small
	enough so that Theorem \ref{A} holds true.

Let us now set
		\begin{equation}
		    b_{\ep}(t) := \m{\uep(t)}, \hspace{2em} 
		    b_{\ep}(0)= |A^{1/2}u_{0}|^{2\gamma}=: b_{0}.
		    \label{defb}
		\end{equation}
 An immediate consequence of Theorem \ref{A} is that
 \begin{equation}
     \frac{K_{3}}{1+t}\leq b_{\ep}(t)\leq \frac{K_{4}}{1+t},
     \hspace{2em} \frac{|b_{\ep}'(t)|}{b_{\ep}(t)}\leq
     \frac{K_{4}}{1+t}\quad\quad\forall t\geq 0,
     \label{h2}
 \end{equation}
 where the constants $K_{3}$ and $K_{4}$ do not depend on $\ep$.
 
Moreover let us define

\begin{equation}
    B_{\ep}(t)=\int_{0}^{t}b_{\ep}(s)\, ds.
    \label{defB}
\end{equation}
From (\ref{h2}) we know that $B_{\ep}(t)\rightarrow +\infty$ as
$t\rightarrow +\infty$, and we have also a non optimal but $\ep$
independent estimate on the speed wherewith  it diverges.  We propose to
obtain sharp estimates on the behaviour of $B_{\ep}$ and, as a 
consequence, also on $\uep$.

 Before proceeding, we introduce some general notations.  Let
 $(e_{k})_{k}$ be a countable basis of $H$ made by eigenvectors of $A$,
 and $\lambda_{k}^{2}$ be the corresponding eigenvalues, that is
 $$Ae_{k} = \lambda_{k}^{2}e_{k},\hspace{2em} 
\forall k;$$
therefore for every $u\in H$ we have that $$u = \sum_{k}u_{k}e_{k}.$$
In particular if $\uep$ is the solution of (\ref{pbm:h-eq}), 
(\ref{pbm:h-data}), then 
$$\uep(t)= \sum_{k}u_{\ep,k}(t)e_{k},$$ 
where $u_{\ep,k}$ solves
\begin{equation}
    \ep u_{\ep,k}''(t) + b_{\ep}(t) \lambda_{k}^{2}u_{\ep,k}(t) + 
    u_{\ep,k}' (t)
    = 0, \hspace{2em} u_{\ep,k}(0) = u_{0,k}, \; u_{\ep,k}'(0) = 
    u_{1,k}. \label{cpcomp} \end{equation}
 Let us now define for $\lambda > 0$:
	%\begin{equation}
 $$
	    H_{\lambda}:=\left\{u\in H: u=\sum_{k:\,\lambda_{k}\geq
	    \lambda}u_{k}e_{k}\right\}, \hspace{2em}
	   H_{\{\lambda\}}:=\left\{u\in H: u=\sum_{k:\,\lambda_{k}=
	    \lambda}u_{k}e_{k}\right\},
	    $$
	%% 
	 %     \label{defhlambda} 
	 % \end{equation}
	 %%
	%\begin{equation}
	$$  H_{[\lambda,\mu)}:=\left\{u\in H: u=\sum_{k:\,\lambda\leq
	    \lambda_{k}< \mu}u_{k}e_{k}\right\}, $$
	%% 
	 %     \label{defhlambda1} 
	 % \end{equation}
	 %%
	and
%	\begin{equation}
	   $$ A_{\lambda}= A_{|H_{\lambda}},\hspace{2em}
	     A_{\{\lambda\}}= A_{|H_{\{\lambda\}}}, \hspace{2em}
	    A_{[\lambda,\mu)}= A_{|H_{[\lambda,\mu)}}. $$
	   % \label{defAlambda} 
	%\end{equation}
Let moreover $\nu$ be defined by:
$$\nu:=\min\{\lambda_{k}: \; u_{0,k}\neq 0, \mbox{ or } 
u_{1,k}\neq 0\}.$$
Since the components of $\uep$ solves (\ref{cpcomp})
then we can assume sine loss of generality that
$\nu^{2}$ is the smallest eigenvalue of A, that is $A = A_{\nu}$, and
that
\begin{equation}
    \langle Au,u\rangle\geq \nu^{2}|u|^{2} \hspace{2em} \forall u\in
    D(A).
    \label{coerc}
\end{equation}
Moreover for every $\mu > \nu$ we can decompose $u\in H$ as 
\begin{equation}
    u= u_{\nu}+ \overline{u}_{\mu} + U_{ \mu}\label{dec}
\end{equation}
where
$u_{\nu}\in H_{\{\nu\}}$, and $ U_{\mu}\in H_{\mu}$.

Finally for every $\lambda \geq \nu$ let us define the corrector
$\Theta_{\ep,\lambda} \in H_{\lambda}$ as the solution of
\begin{equation}
    \ep\Theta_{\ep,\lambda}'' + \Theta_{\ep,\lambda}' = 0,
    \hspace{2em} \Theta_{\ep,\lambda}(0) = 0, \hspace{1em}
    \Theta_{\ep,\lambda}'(0) = U_{1,\lambda}+
    b_{0}A_{\lambda}U_{0,\lambda}.
    \label{pbthetal}
\end{equation}
\subsection{Statements}
We are now ready to state our results.  The first one concerns the
decay of the components of $\uep$.
\begin{thm}
    \label{t1}
    Let $\uep$ be the solution of (\ref{pbm:h-eq}), (\ref{pbm:h-data})
    as in Theorem \ref{A} with $(u_{0},u_{1}) \in D(A)\times D(A^{1/2})$ and let $\lambda \geq \nu$. 
    Then for $\ep$ small (depending
    on $\lambda$) we have the following inequalities.
    \begin{enumerate} 
	\item For $h=0, 1$ there exist constants $\gamma_{h,\lambda}$ 
	independent of $\ep$ such that:
	\begin{equation}
	    \eB{\lambda^{2}}\left(\ep\frac{|A^{h/2}U_{\ep,\lambda}'(t)|^{2}}{b_{\ep}(t)}+
	    |A^{(h+1)/2}U_{\ep,\lambda}(t)|^{2}\right)\leq 
	    \gamma_{h,\lambda}\hspace{2em}\forall t\geq 0.
	    \label{D1}
	\end{equation}
    
        \item  There exists a constant $\gamma_{\lambda}$ 
	independent of $\ep$ such that:
	\begin{equation}
	      \eB{\lambda^{2}}\frac{|U_{\ep,\lambda}'(t)|^{2}}{b_{\ep}^{2}(t)}
	      \leq \gamma_{\lambda} \hspace{2em}\forall t\geq 0.
	    \label{D2}
	\end{equation}
    
        \item  There exists a constant $\gamma_{\ep,\lambda}$ 
	such that:
	\begin{equation}
	      \eB{\lambda^{2}}\frac{|U_{\ep,\lambda}''(t) - 
	      \Theta_{\ep,\lambda}''(t)|^{2}}{b_{\ep}^{4}(t)} \leq 
	      \gamma_{\ep,\lambda} \hspace{2em}\forall t\geq 0.
	    \label{D3}
	\end{equation}
Moreover if $(u_{0},u_{1})\in D(A^{2})\times D(A^{3/2})$, then we can 
take also $\gamma_{\ep,\lambda}$ independent of $\ep$.
    \end{enumerate}
    
\end{thm}

\begin{rmk}
\em{Theorem \ref{t1} is in fact a linear result, indeed in the proof 
we use only that $u_{\ep,k}$ verifies (\ref{cpcomp}) for every $k$ 
with a coefficient $b_{\ep}$ that satisfies (\ref{h2}). This means that 
if the initial data are more regular then estimates like (\ref{D1}) 
- (\ref{D2}) - (\ref{D3}) hold true also for $A^{h/2}U_{\ep,\lambda}$ 
with suitable large $h$.}    
\end{rmk}

Theorem \ref{t1} says that the components of $\uep$ related to big 
eigenvalues decay faster of the component related to the smallest one and as more faster it depends on 
the behaviour of $B_{\ep}$. The following result clarifies this aspect.
\begin{thm}\label{t2}
   Let $\uep$ be the solution of (\ref{pbm:h-eq}), (\ref{pbm:h-data})
    as in Theorem \ref{A} with $(u_{0},u_{1}) \in D(A)\times D(A^{1/2})$. 
    Then for $\ep$ small there exists a non zero vector   $u_{\ep,\infty}\in H_{\{\nu\}}$ 
     such that  as $t\rightarrow +\infty$:
   \begin{equation}
        (1+t)^{1/(2\gamma)}(\uep(t),(1+t)\uep'(t)) \rightarrow (u_{\ep,\infty}, 
    -\frac{1}{2\gamma}u_{\ep,\infty}), \hspace{1em} \mbox{in} \hspace{1em} D(A)\times 
    D(A^{1/2}).
       \label{LIM}
   \end{equation}
   Moreover the following properties hold true.
    \begin{enumerate}
        \item  There exist constants $K_{\ep,1}$, $K_{\ep,2}$ such 
	that for all $t\geq 0$ we have:
	\begin{equation}
	    \frac{K_{\ep,1}}{1+t}\leq e^{-2\nu^{2}\gamma B_{\ep(t)}}\leq  
	    \frac{K_{\ep,2}}{1+t}.
	    \label{B1}
	\end{equation}
   Furthermore if $u_{0,\nu}\neq 0$ then we can take $K_{\ep,1}$ and $K_{\ep,2}$
    independent of $\ep$.
        \item  The following limits hold true for $t\rightarrow 
	+\infty$:
    \begin{equation}
        (1+t)b_{\ep}(t) \rightarrow \frac{1}{2\nu^{2}\gamma};
        \label{B2}
    \end{equation} 
    \begin{equation}
        (1+t)^{1/\gamma}|u_{\ep,\nu}(t)|^{2} \rightarrow 
	\frac{1}{\nu^{2}(2\nu^{2}\gamma)^{1/\gamma}}.
        \label{B31b}
    \end{equation}
    %% 
     % \item The following equalities are verified:
     % \begin{equation}
     %      |u_{\ep,\infty}|^{2} =
     %      \frac{1}{\nu^{2}(2\nu^{2}\gamma)^{1/\gamma}}, \hspace{1em}
     %      |v_{\ep,\infty}|^{2} = 
     %     \frac{\nu^{2}}{(2\nu^{2}\gamma)^{2+1/\gamma}}.
     %     \label{B3b}
     % \end{equation}
     %%
    \item There exists a constant $K_{\ep}$ such that for all $t\geq 
    0$ we have
    \begin{equation}
        |\uep''(t)|^{2}\leq K_{\ep}\frac{1}{(1+t)^{4+1/\gamma}}
        \label{B5}
    \end{equation}
     \end{enumerate}
     \end{thm}
     Let us now give some observations on Theorem \ref{t2}.  
    \begin{itemize}
        \item  Inequalities (\ref{B1}) together  with (\ref{B31b}) and Theorem 
	\ref{t1} say that how much  $U_{\ep,\lambda}$ 
	decay faster of $u_{\ep,\nu}$ depends on $\lambda^{2}/\nu^{2}$.
        \item From limits in (\ref{LIM}) and (\ref{B31b}) we have 
	that  
	$$|u_{\ep,\infty}|^{2} = 
	\frac{1}{\nu^{2}(2\nu^{2}\gamma)^{1/\gamma}},$$
	and it is also obvious that as $t\rightarrow + 
	\infty$ we have that
	\begin{equation}
        (1+t)^{1/\gamma}|A^{1/2}\uep(t)|^{2} \rightarrow 
	\frac{1}{(2\nu^{2}\gamma)^{1/\gamma}}; \hspace{1em}
         (1+t)^{1/\gamma}|A\uep(t)|^{2} \rightarrow 
	\frac{\nu^{2}}{(2\nu^{2}\gamma)^{1/\gamma}};
        \label{B32b}
    \end{equation}
     \begin{equation}
         (1+t)^{2+1/\gamma}|\uep'(t)|^{2} \rightarrow 
	\frac{\nu^{2}}{(2\nu^{2}\gamma)^{2+1/\gamma}}, \hspace{1em} 
	(1+t)^{2+1/\gamma}|A^{1/2}\uep'(t)|^{2} \rightarrow 
	\frac{\nu^{4}}{(2\nu^{2}\gamma)^{2+1/\gamma}},
        \label{B4b}
    \end{equation}
  hence the behaviour at the 
	infinity of the norms do not depends on the initial conditions.
	
    \item   Limits  in (\ref{B32b}) 
	clarify the estimates in (\ref{h1}) and (\ref{h11}) and show that the estimate in 
	(\ref{h12}) is sharp and that a similar estimate holds true 
	(maybe with constants depending on $\ep$) also for 
	$|A^{1/2}\uep'|$. 
    
        \item  The estimate in (\ref{B5}) looks better of all the 
	previous known on the second derivative of $\uep$, moreover 
	it seems optimal, indeed the rate decay is the 
	same as in the limit case $\ep = 0$ with only one real 
	component (in such a case we have only to  solve the ordinary real 
	differential equation $y' + \nu^{2}y^{2\gamma + 1} = 0$).
    
        \item  Theorem \ref{t1} and Theorem \ref{t2} say that as in 
	the case of a linear equation (with a  coefficient $b(t)$ 
	verifying (\ref{h2})) the decay of the solution  is 
	decided only by the smallest eigenvalue $\nu^{2}$ for which there are 
	non zero components of the initial data. Nevertheless in our 
	case the decay rate does not depend on $\nu^{2}$.
    \end{itemize}
	
	\setcounter{equation}{0}
\section{Proofs}\label{sec:proofs}
In some of the proofs we employ the following simple comparison result that
has already been used  in various forms in a lot of 
papers, starting from  \cite{gg:k-dissipative}.

\begin{lemma}\label{lemma:ode}
	Let $f\in C^{1}([0,+\infty))$ and let us
	assume that $f(t)\geq 0$ in $[0,+\infty)$, and that there exist two
	constants $K_{5}>0$, $K_{6}\geq 0$ such that
	$$f'(t)\leq - K_{5}\sqrt{f(t)}\left(\sqrt{f(t)}-K_{6}\right)
	\hspace{2em}
	\forall t\geq 0.$$
	
	Then we have that
	$f(t)\leq\max\left\{f(0),K_{6}^{2}\right\}$ for every
	$t\geq 0$.
\end{lemma}

We divide the proofs in various parts. First we prove two basic 
propositions on linear equations. Then we prove Theorem \ref{t1}. 
After we study the decomposition of $\uep$ made by (\ref{dec}) and 
finally we prove Theorem \ref{t2}.

\subsection{Linear equations and estimates}

Let  $M$ be a self-adjoint linear operator on
$H$. Let us assume that
\begin{equation}
    \langle M w,w\rangle \geq \sigma_{M}^{2}|w|^{2} \hspace{2em} \forall w\in D(M).
    \label{coercM}
\end{equation}
For $h \geq 0$ let us denote  by $|w|_{D(M^{h})}$ the norm of the vector $w$ in 
the space $D(M^{h})$.

Let us assume that $b:[0,+\infty[ \rightarrow ]0,+\infty[$  is a 
$C^{1}$ function that verifies 
 \begin{equation}
    b(t)\leq \frac{K_{4}}{1+t},
     \hspace{2em} \frac{|b'(t)|}{b(t)}\leq
     \frac{K_{4}}{1+t}, \hspace{2em} \frac{|b'(t)|}{b^{2}(t)}\leq
   \frac{ K_{4}}{K_{3}}\quad\quad\forall t\geq 0,
     \label{h2b}
 \end{equation}
 where $K_{4}$ and $K_{3}$ are the constants in (\ref{h2}).
For simplicity in the following we use these notations:
$$\|u\|^{2} = |u|^{2}\left(1 + b(0)^{-1}+b(0)^{-2}\right)\,\, \mbox{ if 
$u\in H$,}$$

 $$\|u\|_{D(M^{h/2})}^{2} = |u|_{D(M^{h/2})}^{2}\left(1 + 
 b(0)^{-1}+b(0)^{-2}\right)\,\, \mbox{ if 
$u\in D(M^{h/2})$.}$$
Let $v_{\ep}  \in C^{2}([0,+\infty[,H) \cap  
C^{1}([0,+\infty[,D(M^{1/2})) \cap  C^{0}([0,+\infty[,D(M))$ be the 
solution of the problem:
\begin{equation}
    \ep v_{\ep}''(t) +b(t) Mv_{\ep}(t) + v_{\ep}'(t) = 0, \hspace{2em} v_{\ep}(0) = 
    v_{0}\in D(M), \hspace{1em} v_{\ep}'(0) = v_{1}\in D(M^{1/2}).
    \label{cpM}
\end{equation}
Moreover let $\theta_{\ep}$ be the solution of
\begin{equation}
    \ep \theta_{\ep}''(t) +\theta_{\ep}'(t) = 0, \hspace{2em} \theta_{\ep}(0) = 
    0, \hspace{1em} \theta_{\ep}'(0) = v_{1} + b(0) Mv_{0},
    \label{cptheta}
\end{equation}
so that $\theta_{\ep}(t) = \ep\theta_{\ep}'(0)(1 - e^{-t/\ep})$, and let us set 
$$w_{\ep} = v_{\ep} - \theta_{\ep}.$$

Finally let $B$ defined as in (\ref{defB}) (using $b(t)$ in place of 
$b_{\ep}$ of course).

Therefore the following propositions hold true.
\begin{prop}
    \label{p1}
     Let $h\geq 1$ and let us assume 
    that $(v_{0},v_{1})\in D(M^{(h+1)/2})\times D(M^{h/2})$. Then for 
    $\ep $ small depending only on $\sigma_{M}^{2}$, $K_{3}$ and 
    $K_{4}$
    (and not on the initial data or $h$), for all $t\geq 0$ we have 
    that:
    \begin{equation}
        e^{2\sigma_{M}^{2}B(t)}
	\left(\ep\frac{|M^{h/2}v_{\ep}'(t)|^{2}}{b(t)}+
	    |M^{(h+1)/2}v_{\ep}(t)|^{2}\right)\leq 
	    L_{0}(\|v_{1}\|_{D(M^{h/2})}^{2} + 
	    |v_{0}|_{D(M^{(h+1)/2})}^{2})=:L_{h,M},
        \label{SL1b}
    \end{equation}
    \begin{equation}
       e^{2\sigma_{M}^{2}B(t)}
	\frac{|M^{(h-1)/2}v_{\ep}'(t)|^{2}}{b^{2}(t)}\leq 
	    L_{1}(\|v_{1}\|_{D(M^{h/2})}^{2} + 
	    |v_{0}|_{D(M^{(h+1)/2})}^{2}) =:H_{h,M},
        \label{SL2b}
    \end{equation}
    where $L_{0}$ and $L_{1}$ depend only on $\sigma_{M}^{2}$, $K_{3}$ and 
    $K_{4}$.
    \end{prop}
    
\begin{prop}
    \label{p2}
    Let  us assume 
    that $(v_{0},v_{1})\in D(M^{2})\times D(M^{3/2})$. Then for 
    $\ep $ small depending only on  $\sigma_{M}^{2}$, $K_{3}$ and 
    $K_{4}$ 
    and not on the initial data we have that
    \begin{equation}
       e^{2\sigma_{M}^{2}B(t)}
\frac{ |w_{\ep}''(t)|^{2}}{b^{4}(t)}\leq 
	    L_{2}(\|v_{1}\|_{D(M^{3/2})}^{2} + 
	    |v_{0}|_{D(M^{2})}^{2}),\hspace{2em} \forall t\geq 0,
        \label{SL3b}
    \end{equation}
    where $L_{2}$ depends only on $\sigma_{M}^{2}$, $K_{3}$ and 
    $K_{4}$.
    \end{prop}
    
    Now let us prove Proposition \ref{p1} and Proposition \ref{p2}.
    \paragraph{Proof of Proposition \ref{p1}}
    Let us denote by $c_{i}$ and $C_{i}$ various constants that depend only on 
     $\sigma_{M}^{2}$, $K_{3}$  and $K_{4}$. 
     
     The outline of the proof is the following. Firstly ({\em Step 
     1}), we prove,  for 
     every $h\geq 0$, that 
     
   \emph{if  we have that
     \begin{equation}
         e^{2\sigma_{M}^{2}B(t)} |M^{h/2}v_{\ep}(t)|^{2}\leq R_{h} 
	 \hspace{2em} \forall t\geq 0
         \label{fond}
     \end{equation}
     then for all $t\geq 0$ we get that
     \begin{equation}
          e^{2\sigma_{M}^{2}B(t)}\left[\ep \frac{|M^{h/2}v_{\ep}'(t)|^{2}}{b(t)} +
        |M^{(h+1)/2}v_{\ep}(t)|^{2}\right]\leq 16\sigma_{M}^{2}R_{h} + 
	C_{0} (\|v_{1}\|_{D(M^{h/2})}^{2} + 
	|v_{0}|_{D(M^{(h+1)/2})}^{2}).
         \label{sl}
     \end{equation}}
     Since problem (\ref{cpM}) is linear it is enough to  prove this 
     estimate for $h= 0$.
     
     Then ({\em Step 2}) we show that for $h=1$ we have (\ref{fond}) with
    \begin{equation}
        R_{1}= C_{1}(\|v_{1}\|_{D(M^{1/2})}^{2} + 
	|v_{0}|_{D(M)}^{2}).
        \label{r1}
    \end{equation}
     Using (\ref{r1}) in (\ref{sl}) with $h=1$ we can then conclude that 
     (\ref{SL1b}) holds true if $h=1$. Since (\ref{cpM}) is linear, 
     (\ref{SL1b}) will proved for every $h\geq 1$.
    
   In conclusion for proving (\ref{SL1b}) we have only to prove 
   (\ref{sl}) with $h=0$ and 
   (\ref{r1}).
   
  Finally ({\em Step 3}) we prove (\ref{SL2b}). Also in this case it is 
   enough to consider the case $h=1$.

   For $\alpha > 0$  let us introduce the following energies that we 
   use in the proofs:
    \begin{eqnarray*}
        D_{\alpha}(t) & := & e^{2\alpha B(t)}\left[\langle \ep v_{\ep}'(t),
       v_{\ep}(t)\rangle +\frac{1 }{2}|v_{\ep}(t)|^{2}\right], \\
        E_{\alpha}(t) & := & e^{2\alpha B(t)}\left[\ep \frac{|v_{\ep}'(t)|^{2}}{b(t)} +
        |M^{1/2}v_{\ep}(t)|^{2}\right] , \\
        F_{\alpha}(t) & :=  &e^{2\alpha B(t)}  \frac{|v_{\ep}'(t)|^{2}}{b^{2}(t)}.
    \end{eqnarray*}
An easy calculation shows that
\begin{equation}
    D_{\alpha}' = 2\alpha b D_{\alpha} -b e^{2\alpha B}|M^{1/2}v_{\ep}|^{2}
   + \ep e^{2\alpha B}|v_{\ep}'|^{2};
    \label{derD}
\end{equation}
\begin{equation}
    E_{\alpha}'=  -e^{2\alpha B} \frac{|v_{\ep}'|^{2}}{b}\left(2 
    +\ep\frac{b'}{b} -2\alpha\ep b\right)+ 2\alpha b e^{2\alpha B}
    |M^{1/2}v_{\ep}|^{2};
    \label{derE}
\end{equation}
\begin{equation}
    F_{\alpha}'= -\frac{1}{\ep}F_{\alpha} \left(2 
    +2\ep\frac{b'}{b} -2\alpha\ep b\right) - \frac{2}{\ep}e^{2\alpha 
    B}\frac{1}{b}\langle v_{\ep}',Mv_{\ep}\rangle.
    \label{derF}
\end{equation}

\subparagraph{\emph{Step 1 - Proof of (\ref{sl}) with $h=0$}}

Let us choose $\alpha = 2\sigma_{M}^{2}:=\alpha_{0}$.

\subparagraph{\textsc{Estimate on $D_{\alpha_{0}}$}}
We prove that, if $\ep$ is small enough, for all $t\geq 0$ we have 
that
\begin{eqnarray}
    \int_{0}^{t}e^{2\alpha_{0} B(s)}b(s) |M^{1/2}v_{\ep}(s)|^{2}ds &\leq &
   |v_{1}|^{2}+|v_{0}|^{2} + 
  C_{2} \ep^{2}e^{2\alpha_{0} B(t)}\frac{|v_{\ep}'(t)|^{2}}{b(t)} + \nonumber \\
    &&
   + C_{3}\ep\int_{0}^{t} e^{2\alpha_{0} B(s)}  
    \frac{|v_{\ep}'(s)|^{2}}{b(s)} ds + 2R_{0}e^{\alpha_{0} B(t)}.
    \label{stD3-1}
\end{eqnarray}
By (\ref{fond}) we obtain that 
\begin{eqnarray}
    2\alpha_{0} b D_{\alpha_{0}} & = & 2\alpha_{0}\ep b e^{2\alpha_{0} B}\langle 
    v_{\ep}',v_{\ep}\rangle +\alpha_{0} b e^{2\alpha_{0} B}|v_{\ep}|^{2}
    \nonumber \\
     &  \leq & \alpha_{0}\ep^{2} b e^{2\alpha_{0} B}|v_{\ep}'|^{2}+2\alpha_{0} b 
     e^{2\alpha_{0} B}|v_{\ep}|^{2}
    \nonumber  \\
     & \leq & \alpha_{0} \ep^{2} b^{2}e^{2\alpha_{0} B}  
    \frac{|v_{\ep}'|^{2}}{b} + 2\alpha_{0} R_{0}b 
     e^{\alpha_{0} B}.
    \label{stD0-1}
\end{eqnarray}
From (\ref{derD}) and (\ref{stD0-1}) we therefore get that
\begin{equation}
    D_{\alpha_{0}}' + 
     e^{2\alpha_{0} B}b|M^{1/2}v_{\ep}|^{2}\leq 
   \ep (b + \alpha_{0} \ep  b^{2})e^{2\alpha_{0} B}  
    \frac{|v_{\ep}'|^{2}}{b} + 2\alpha_{0} R_{0}b 
     e^{\alpha_{0} B}. \label{D-1-1}
\end{equation}
Since from (\ref{h2b}) the function $b$ is bounded by $K_{4}$ then integrating 
(\ref{D-1-1}) we arrive at
\begin{equation}
    \int_{0}^{t}
     e^{2\alpha_{0} B(s)}b(s)|M^{1/2}v_{\ep}(s)|^{2}ds \leq D_{\alpha_{0}}(0) 
     - D_{\alpha_{0}}(t) + 
     \ep c_{1}\int_{0}^{t} e^{2\alpha_{0} B(s)}  
    \frac{|v_{\ep}'(s)|^{2}}{b(s)} ds+ 2R_{0}e^{\alpha_{0} B(t)}.
    \label{stD-1-1}
\end{equation}
Since $\ep\leq 1$ and $b$ is bounded, we can estimate $D_{\alpha_{0}}(0)$ and $D_{\alpha_{0}}(t)$ as follows:
\begin{eqnarray}
    |D_{\alpha_{0}}(0)| & \leq & \ep|v_{1}||v_{0}| + \frac{1}{2}|v_{0}|^{2} \leq 
    |v_{1}|^{2}+|v_{0}|^{2},
    \label{dati1}  \\
    -D_{\alpha_{0}}(t) & \leq  & e^{2\alpha_{0} B(t)}\left (\ep |v_{\ep}'(t)||v_{\ep}(t)| -  
    \frac{1}{2}|v_{\ep}(t)|^{2}\right) \nonumber
    \\&\leq &
    \frac{1}{2}\ep^{2}e^{2\alpha_{0} 
    B(t)}\frac{|v_{\ep}'(t)|^{2}}{b(t)}b(t)\leq  c_{2}\ep^{2}e^{2\alpha_{0} 
    B(t)}\frac{|v_{\ep}'(t)|^{2}}{b(t)}.
    \label{stD-2-1}
\end{eqnarray}
Plugging (\ref{dati1}) - (\ref{stD-2-1}) in (\ref{stD-1-1}) 
 we achieve (\ref{stD3-1}).

 \subparagraph{\textsc{Proof of (\ref{sl})}}
Integrating (\ref{derE}) and using (\ref{stD3-1}) we get that
\begin{eqnarray}
    E_{\alpha_{0}}(t) & \leq & E_{\alpha_{0}}(0) - \int_{0}^{t}e^{2\alpha_{0} B(s)}  
    \frac{|v_{\ep}'(s)|^{2}}{b(s)} \left(2 +\ep 
    \frac{b'(s)}{b(s)}-2\alpha_{0}\ep b(s) -2\alpha_{0} C_{3}\ep \right) ds+
    \nonumber  \\
     &  & +2\alpha_{0} C_{2}\ep^{2}e^{2\alpha_{0} B(t)}
     \frac{|v_{\ep}'(t)|^{2}}{b(t)}
     + 2\alpha_{0} ( |v_{1}|^{2}+|v_{0}|^{2})+ 4\alpha_{0} R_{0}e^{\alpha_{0} B(t)}.
    \label{E-1-1}
\end{eqnarray}
Thanks to (\ref{h2b}) we can take $\ep$ small enough in such a way 
that 
\begin{equation}
    2 - 2 \ep\sup_{t\geq 0}\frac{|b'(t)|}{b(t)} 
    -4\sigma_{M}^{2}\ep\sup_{t\geq 0} b(t) 
    -4\sigma_{M}^{2}\ep C_{3} \geq 1,
    \label{ep21}
\end{equation}
\begin{equation}
   4\sigma_{M}^{2} C_{2}\ep \leq \frac{1}{2}.
    \label{ep31}
\end{equation}
Plugging (\ref{ep21}) and (\ref{ep31}) in (\ref{E-1-1}) we obtain that
$$E_{\alpha_{0}} (t) \leq \frac{|v_{1}|^{2}}{b(0)} + |M^{1/2}v_{0}|^{2}
+ \frac{1}{2}E_{\alpha_{0}}(t)+ 
2\alpha_{0}(|v_{1}|^{2}+|v_{0}|^{2})  +4\alpha_{0} R_{0}e^{\alpha_{0} B(t)}$$
from which 
$$ \frac{1}{2} E_{\alpha_{0}}(t) \leq c_{3}
(\|v_{1}\|^{2}+|v_{0}|_{D(M^{1/2})}^{2})  +4\alpha_{0} 
R_{0}e^{\alpha_{0} B(t)}.$$
Since $\alpha_{0} = 2\sigma_{M}^{2}$, hence (\ref{sl}) immediately follows dividing all terms by 
$e^{\alpha_{0} B(t)}$.    

\subparagraph{\emph{Step 2 - Proof of (\ref{r1})}}
For begin with, let us choose 
\begin{equation}
  \alpha = \sigma_{M}^{2} - \frac{1}{8K_{4}}:= \beta. 
   \label{defbeta}
\end{equation}

Firstly we prove that for $\ep$ small and $h= 0$, $h=1$ we have for all 
$t\geq 0$ that:
\begin{equation}       
    e^{2\beta B(t)}
	\left(\ep\frac{|M^{h/2}v_{\ep}'(t)|^{2}}{b(t)}+
	    |M^{(h+1)/2}v_{\ep}(t)|^{2}\right)\leq 
	    C_{4}(\|v_{1}\|_{D(M^{h/2})}^{2} + 
	    |v_{0}|_{D(M^{(h+1)/2})}^{2}),
        \label{SL1}
    \end{equation}
    and 
    \begin{equation}
         e^{2\beta B(t)}
	\frac{|v_{\ep}'(t)|^{2}}{b^{2}(t)}\leq 
	    C_{5}(\|v_{1}\|_{D(M^{1/2})}^{2} + 
	    |v_{0}|_{D(M)}^{2}).
        \label{SL2}
    \end{equation}
Since (\ref{cpM}) is linear it is enough to prove (\ref{SL1}) with 
$h=0$.

\subparagraph{\textsc{Estimate on $D_{\beta}$}}
We prove that, if $\ep$ si small enough, for all $t\geq 0$ we have 
that
\begin{eqnarray}
    \int_{0}^{t}e^{2\beta B(s)}b(s) |M^{1/2}v_{\ep}(s)|^{2}ds &\leq &
    C_{6}(|v_{1}|^{2}+|v_{0}|^{2}) + 
    C_{7}\ep^{2}e^{2\beta B(t)}\frac{|v_{\ep}'(t)|^{2}}{b(t)} + \nonumber \\
    &&
   + C_{8}\ep\int_{0}^{t} e^{2\beta B(s)}  
    \frac{|v_{\ep}'(s)|^{2}}{b(s)} ds.
    \label{stD3}
\end{eqnarray}
From (\ref{coercM}) we obtain that 
\begin{eqnarray}
    2\beta b D_{\beta} & = & 2\beta\ep b e^{2\beta B}\langle 
    v_{\ep}',v_{\ep}\rangle +\beta b e^{2\beta B}|v_{\ep}|^{2}
    \nonumber \\
     &  \leq & \beta\ep b e^{2\beta B}|v_{\ep}'|^{2}+\beta(1+\ep)b 
     e^{2\beta B}|v_{\ep}|^{2}
    \nonumber  \\
     & \leq & \beta \ep b^{2}e^{2\beta B}  
    \frac{|v_{\ep}'|^{2}}{b} + \frac{\beta}{\sigma_{M}^{2}}(1+\ep)b 
     e^{2\beta B}|M^{1/2}v_{\ep}|^{2}.
    \label{stD0}
\end{eqnarray}
From (\ref{derD}) and (\ref{stD0}) we therefore get that
\begin{equation}
    D_{\beta}' + \left(1 - \frac{\beta}{\sigma_{M}^{2}}(1+\ep)\right)
     e^{2\beta B}b|M^{1/2}v_{\ep}|^{2}\leq 
   \ep (b + \beta b^{2})e^{2\beta B}  
    \frac{|v_{\ep}'|^{2}}{b}. \label{D-1}
\end{equation}
Since by (\ref{h2b}) the function $b$ is bounded then integrating 
(\ref{D-1}) we arrive at
\begin{eqnarray}
    &\displaystyle \left(1 - \frac{\beta}{\sigma_{M}^{2}}(1+\ep)\right)\int_{0}^{t}
     e^{2\beta B(s)}b(s)|M^{1/2}v_{\ep}(s)|^{2}ds \leq D_{\beta}(0) - 
     D_{\beta}(t) + &\nonumber \\
     &
     +\ep c_{4}\displaystyle\int_{0}^{t} e^{2\beta B(s)}  
    \frac{|v_{\ep}'(s)|^{2}}{b(s)} ds.&
    \label{stD-1}
\end{eqnarray}
We can estimate $D_{\beta}(0)$
and $D_{\beta}(t)$ as  is (\ref{dati1})  and   (\ref{stD-2-1}),
furthermore since $\beta  < \sigma_{M}^{2}$ we can take $\ep$ small 
enough in such a way that
\begin{equation}
    1 - \frac{\beta }{\sigma_{M}^{2}}(1+\ep) \geq c_{5} > 0.
    \label{ep1}
\end{equation}
Plugging (\ref{ep1}), (\ref{dati1}), (\ref{stD-2-1})
(with $\beta$ instead of $\alpha_{0}$) in (\ref{stD-1}) 
 we achieve (\ref{stD3}).

 \subparagraph{\textsc{Proof of (\ref{SL1}) with $h=0$}}
Integrating (\ref{derE}) and using (\ref{stD3}) we get that
\begin{eqnarray}
    E_{\beta }(t) & \leq & E_{\beta}(0) - \int_{0}^{t}e^{2\beta  B(s)}  
    \frac{|v_{\ep}'(s)|^{2}}{b(s)} \left(2 +\ep 
    \frac{b'(s)}{b(s)}-2\beta\ep b(s) -2\beta \ep C_{8}\right) ds+
    \nonumber  \\
     &  & +2\beta  C_{7}\ep^{2}e^{2\beta  B(t)}
     \frac{|v_{\ep}'(t)|^{2}}{b(t)}
     + 2\beta  C_{6} ( |v_{1}|^{2}+|v_{0}|^{2}).
    \label{E-1}
\end{eqnarray}
Thanks to (\ref{h2b}) we can take $\ep$ small enough in such a way 
that 
\begin{equation}
    2 - 2 \ep\sup_{t\geq 0}\frac{|b'(t)|}{b(t)} -2\beta 
    \ep\sup_{t\geq 0} b(t) 
    -2\beta \ep C_{8} \geq 1,
    \label{ep2}
\end{equation}
\begin{equation}
    2\beta  C_{7}\ep \leq \frac{1}{2}.
    \label{ep3}
\end{equation}
Plugging (\ref{ep2}) and (\ref{ep3}) in (\ref{E-1}) we obtain that
$$E_{\beta }(t) \leq \frac{|v_{1}|^{2}}{b(0)} + |M^{1/2}v_{0}|^{2} + 
2\beta C_{6} (|v_{1}|^{2}+|v_{0}|^{2}) + \frac{1}{2}E_{\beta}(t),$$
from which (\ref{SL1}) immediately follows.    

\subparagraph{\textsc{Proof of (\ref{SL2})}}
Plugging (\ref{ep2})  in (\ref{derF}) we have that
$$F_{\beta }' \leq -\frac{1}{\ep} F_{\beta } +
\frac{2}{\ep}\sqrt{F_{\beta}}|Mv_{\ep}|e^{\beta  B}.$$
Applying (\ref{SL1}) with $h=1$ we then obtain that
$$F_{\beta}' \leq -\frac{1}{\ep} \sqrt{F_{\beta}}\left(\sqrt{F_{\beta }} - 
2\sqrt{C_{4}(\|v_{1}\|_{D(M^{1/2})}^{2} + 
	    |v_{0}|_{D(M)}^{2})}\right),$$
hence from Lemma \ref{lemma:ode} we get that
$$F_{\beta}(t) \leq \max\{F_{\beta}(0), 4C_{4}(\|v_{1}\|_{D(M^{1/2})}^{2} + 
	    |v_{0}|_{D(M)}^{2})\}, \quad\quad \forall t \geq 0,$$
that is (\ref{SL2}).

\subparagraph{Proof of (\ref{r1})} Since (\ref{sl}) holds true, it is 
enough to prove that (\ref{fond}) holds true with $h=0$ and 
\begin{equation}
         R_{0} =  C_{9}(\|v_{1}\|_{D(M^{1/2})}^{2} + 
	    |v_{0}|_{D(M)}^{2}).
    \label{r0}
\end{equation}
To this end let us set $\alpha = \sigma_{M}^{2}$. 
By (\ref{coercM}) we obtain that 
\begin{eqnarray}
    2\sigma_{M}^{2} b D_{\sigma_{M}^{2}} & = & 2\sigma_{M}^{2}
    \ep b e^{2\sigma_{M}^{2} B}\langle 
    v_{\ep}',v_{\ep}\rangle +\sigma_{M}^{2} b e^{2\sigma_{M}^{2} B}|v_{\ep}|^{2}
    \nonumber \\
     &  \leq & \sigma_{M}^{2}\ep^{2}  e^{2\sigma_{M}^{2} 
     B}|v_{\ep}'|^{2}+\sigma_{M}^{2}b^{2} 
     e^{2\sigma_{M}^{2} B}|v_{\ep}|^{2} + \sigma_{M}^{2}b
     e^{2\sigma_{M}^{2} B}|v_{\ep}|^{2}
    \nonumber  \\
     & \leq &\sigma_{M}^{2}\ep^{2}  e^{2\sigma_{M}^{2}
     B}|v_{\ep}'|^{2} + b^{2} 
     e^{2\sigma_{M}^{2} B}|M^{1/2}v_{\ep}|^{2} +b 
     e^{2\sigma_{M}^{2} B}|M^{1/2}v_{\ep}|^{2}.
    \label{stD00}
\end{eqnarray}
Moreover from (\ref{h2b})  and (\ref{defbeta}) we have that
\begin{equation}
    e^{2\sigma_{M}^{2}B(t)}= e^{2\beta B(t)}e^{2(\sigma_{M}^{2 }- 
    \beta)B(t)} \leq e^{2\beta B(t)} e^{2K_{4}(\sigma_{M}^{2}- 
    \beta)\log(1+t)} = e^{2\beta B(t)} (1+t)^{1/4},
    \label{exps}
\end{equation}
hence using once again (\ref{h2b}) and (\ref{SL1}) with $h=0$, 
inequality
(\ref{stD00}) becomes

\begin{eqnarray}
     2\sigma_{M}^{2} b D_{\sigma_{M}^{2}}&\leq& 
     \sigma_{M}^{2}\ep^{2}  e^{2\sigma_{M}^{2}
     B}|v_{\ep}'|^{2} + c_{6}(1+t)^{-7/4} 
     e^{2\beta B}|M^{1/2}v_{\ep}|^{2} +b 
     e^{2\sigma_{M}^{2} B}|M^{1/2}v_{\ep}|^{2}\nonumber
    \\
    &\leq &\sigma_{M}^{2}\ep^{2}  e^{2\sigma_{M}^{2}
     B}|v_{\ep}'|^{2} + c_{7}(1+t)^{-7/4} 
     (\|v_{1}\|^{2}+ |v_{0}|_{D(M^{1/2})}^{2}) +\nonumber \\
     & & + b 
     e^{2\sigma_{M}^{2} B}|M^{1/2}v_{\ep}|^{2}.
      \label{MD-1}
\end{eqnarray}

Plugging (\ref{MD-1}) into (\ref{derD}) and integrating we obtain that
\begin{equation}    D_{\sigma_{M}^{2}}(t)\leq D_{\sigma_{M}^{2}}(0)+ 
    c_{8} 
     (\|v_{1}\|^{2}+ |v_{0}|_{D(M^{1/2})}^{2}) +\ep\int_{0}^{t}
e^{2\sigma_{M}^{2}
     B(s)}\frac{|v_{\ep}'(s)|^{2}}{b^{2}(s)}b^{2}(s) (1 + 
     \ep\sigma_{M}^{2})\, ds  .
         \label{Din}
     \end{equation}
     From (\ref{exps}), (\ref{SL2}) and (\ref{h2b})  we get that
     \begin{eqnarray}
         e^{2\sigma_{M}^{2}
     B(t)}\frac{|v_{\ep}'(t)|^{2}}{b^{2}(t)}b^{2}(t)&\leq &
     e^{2\beta
     B(t)}\frac{|v_{\ep}'(t)|^{2}}{b^{2}(t)}b^{2}(t)(1+t)^{1/4}\nonumber\\
     &\leq & C_{5}(\|v_{1}\|_{D(M^{1/2})}^{2} + 
     |v_{0}|_{D(M)}^{2})b^{2}(t)(1+t)^{1/4} \nonumber \\
     &\leq& c_{9}(\|v_{1}\|_{D(M^{1/2})}^{2} + 
     |v_{0}|_{D(M)}^{2})(1+t)^{-7/4}.
         \label{stv'}
     \end{eqnarray}
Plugging (\ref{stv'}) into  (\ref{Din}) we arrive at
\begin{eqnarray*}
    \frac{1}{2} e^{2\sigma_{M}^{2} B(t)}|v_{\ep}(t)|^{2}& \leq & 
    |D_{\sigma_{M}^{2}}(0)| +
    \ep e^{2\sigma_{M}^{2} B} |\langle v_{\ep}'(t), v_{\ep}(t) \rangle|
    + c_{8} 
     (\|v_{1}\|^{2}+ |v_{0}|_{D(M^{1/2})}^{2})  +  \\
      & & +\ep c_{10} 
     (\|v_{1}\|_{D(M^{1/2})}^{2}+ |v_{0}|_{D(M)}^{2})  \\
     & \leq & c_{11} (\|v_{1}\|_{D(M^{1/2})}^{2}+ |v_{0}|_{D(M)}^{2})+
     \ep^{2}e^{2\sigma_{M}^{2}
     B(t)}\frac{|v_{\ep}'(t)|^{2}}{b^{2}(t)}b^{2}(t) + \frac{1}{4}
     e^{2\sigma_{M}^{2}
     B(t)}|v_{\ep}(t)|^{2}  \\
     & \leq & c_{12}(\|v_{1}\|_{D(M^{1/2})}^{2}+ |v_{0}|_{D(M)}^{2}) + \frac{1}{4}
     e^{2\sigma_{M}^{2}
     B(t)}|v_{\ep}(t)|^{2}.
\end{eqnarray*}
By this last inequality (\ref{r0}) immediately follows.

\subparagraph{\emph{Step 3 - Proof of (\ref{SL2b}) with $h=1$}}
Let $\alpha = \sigma_{M}^{2}$. From (\ref{derF}), (\ref{ep21}) and 
(\ref{SL1b}) used with $h=1$ we deduce that
$$F_{\sigma_{M}^{2}}' \leq -\frac{1}{\ep}\sqrt{F_{\sigma_{M}^{2}}}
\left(\sqrt{F_{\sigma_{M}^{2}}} - 2|Mv_{\ep}|e^{\sigma_{M}^{2}B}\right) \leq 
-\frac{1}{\ep}\sqrt{F_{\sigma_{M}^{2}}}\left (\sqrt{F_{\sigma_{M}^{2}}} - 
2\sqrt{L_{1,M}}\right).$$
We can then apply Lemma \ref{lemma:ode}, hence for all $t\geq 0$ we 
have that
$$F_{\sigma_{M}^{2}}(t) \leq \max \{ F_{\sigma_{M}^{2}}(0), 4L_{1,M}\}
\leq F_{\sigma_{M}^{2}}(0) + 4L_{1,M},$$
therefore (\ref{SL2b}) holds true.

\qed
\paragraph{Proof of Proposition \ref{p2}}
   Let us take $\ep$ small enough in 
   such a way that we can apply Proposition \ref{p1} (with $h= 3$ and 
   $h=1$).
  
   Firstly let us observe that  $w_{\ep}$ satisfies the following 
   problem
 \begin{equation}
     \ep w_{\ep}''(t) + w_{\ep}'(t) = -b(t)Mv_{\ep}(t), \hspace{1em}
     w_{\ep}(0) = v_{0}, \hspace{0,5em} w_{\ep}'(0) = -b(0)Mv_{0}, 
     \hspace{0,5em} w_{\ep}''(0) = 0.
     \label{eqw}
 \end{equation}
If we  set
 $$G(t):= e^{2\sigma_{M}^{2}B(t)}\frac{|w_{\ep}''(t)|^{2}}{b^{4}(t)},$$
therefore from (\ref{eqw}) we have that
$$G'= G\left(2\sigma_{M}^{2}b - 4\frac{b'}{b}\right) - \frac{2}{\ep}
\frac{e^{2\sigma_{M}^{2}B}}{b^{4}}\langle w_{\ep}'', -w_{\ep}'' 
-bMv_{\ep}' -b'Mv_{\ep}\rangle.$$
Hence we immediately get that
\begin{equation}
    G' \leq -\frac{1}{\ep}G\left(2- 2\sigma_{M}^{2}b\ep + 4\ep
   \frac{b'}{b}\right) + \frac{2}{\ep} 
   \sqrt{G}\left(\frac{|Mv_{\ep}'|}{b} + \frac{|b'|}{b^{2}} 
   |Mv_{\ep}|\right)e^{\sigma_{M}^{2}B}.
    \label{derG}
\end{equation}
Thanks to (\ref{h2b}) we can take $\ep$ small enough so that
\begin{equation}
    2 -  2 \ep\sigma_{M}^{2}\sup_{t\geq 0}b(t) - 4\ep
   \sup_{t\geq 0}\frac{|b'(t)|}{b(t)}\geq 1.
    \label{ep7}
\end{equation}
Using (\ref{ep7}), (\ref{h2b}), (\ref{SL1b}) with $h=1$ and 
   (\ref{SL2b}) with $h=3$ in (\ref{derG})
   we obtain that
   $$G' \leq -\frac{1}{\ep} \sqrt{G}\left( \sqrt{G} - 
   c_{1}(\|v_{1}\|_{D(M^{3/2})}^{2} + 
   |v_{0}|_{D(M^{2})}^{2})^{1/2}\right),$$
   with a constant $c_{1}$ that depends only on 
   $\sigma_{M}^{2}$ and $K_{3}$, $K_{4}$.
  Thus we can apply Lemma \ref{lemma:ode}, from which we have that
   $$G(t) \leq \max\{G(0),  c_{1}^{2}(\|v_{1}\|_{D(M^{3/2})}^{2} + 
   |v_{0}|_{D(M^{2})}^{2})\}, \quad \quad \forall t\geq 0.$$
   Since $G(0) = 0$ thesis is proved. 
   \qed
   \subsection{Proof of Theorem \ref{t1}}
  We denote by $c_{i,\lambda}$ various constants that depend only on 
  $\lambda$ and   on $|u_{0}|_{D(A)}$, $|u_{1}|_{D(A^{1/2})}$.
  
  To begin with let us recall that thanks to (\ref{h2})  functions 
  $b_{\ep}$ verify (\ref{h2b}) independently of $\ep$. Let us also 
  stress that by (\ref{defb}) we have $b_{\ep}(0) = b_{0}$ 
  independent of $\ep$.
  
  To obtain inequalities (\ref{D1})  and (\ref{D2}) it is enough to 
   apply Proposition \ref{p1} with $M=A_{\lambda}$, $b(t) = 
   b_{\ep}(t)$ and $\sigma_{M}^{2} = \lambda^{2}$
   (taking of course $\ep$ small enough); indeed in such a case 
   $U_{\ep,\lambda}$ solves (\ref{cpM}).
   
   Now let us prove (\ref{D3}). 
   
   When the initial data are regular we 
   can apply directly Proposition \ref{p2} with $M= A_{\lambda}$ and 
   we obtain (\ref{D3}) with a constant that does not depend on $\ep$.
   
   Now 
   let us consider the general case in which $(u_{0},u_{1})\in 
   D(A)\times D(A^{1/2})$. Let us set
   $$\mu^{2} := \lambda^{2} + \frac{1}{K_{3}}$$
   where $K_{3}$ is the constant in (\ref{h2}). Then we can write
   $$U_{\ep,\lambda}= V_{\ep,\lambda} + U_{\ep,\mu}, \hspace{2em}
   \Theta_{\ep,\lambda}= \theta_{\ep,\lambda} + \Theta_{\ep,\mu}.$$
   We estimate separately $V_{\ep,\lambda}''-  \theta_{\ep,\lambda}''$ 
   and $U_{\ep,\mu}''- \Theta_{\ep,\mu}''$.
   \paragraph{\textsc{Estimate on $U_{\ep,\mu}''- \Theta_{\ep,\mu}''$}}
   We prove that for every $t\geq 0$ we have that:
   \begin{equation}
       e^{2\lambda^{2}B_{\ep}(t)}\frac{1}{b_{\ep}^{4}(t)}|U_{\ep,\mu}''(t) - 
       \Theta_{\ep,\mu}''(t)|^{2}\leq 
       \frac{c_{1,\lambda}}{\ep^{2}}.
       \label{S3}
   \end{equation}
Let us assume that $M= A_{\mu}$, $b(t) = b_{\ep}(t)$ and that $\ep$ 
is small enough in such a way that we can apply Proposition \ref{p1} 
with these choices. Then from (\ref{cpM}), (\ref{SL1b}), (\ref{SL2b}) with $h=1$ 
and (\ref{h2}) we obtain that:
\begin{eqnarray}
    \eB{\lambda^{2}}|U_{\ep,\mu}''(t)|^{2} &\leq &\frac{2}{\ep^{2}} 
    \eB{(\lambda^{2}-\mu^{2})}  \eB{\mu^{2}}(b_{\ep}^{2}(t)|M 
    U_{\ep,\mu}(t)|^{2} + |U_{\ep,\mu}'(t)|^{2}) \nonumber \\
    &\leq & \frac{2b_{\ep}^{2}(t)}{\ep^{2}} 
    (L_{0}+L_{1})(\|u_{1}\|_{D(A^{1/2})}^{2}+|u_{0}|_{D(A)}^{2}) 
    e^{2(\lambda^{2}-\mu^{2})K_{3}\log(1+t)} \nonumber \\
    &= & \frac{1}{\ep^{2}}c_{2,\lambda}b_{\ep}^{2}(t)(1+t)^{-2}.
    \label{umu2}
\end{eqnarray}
Moreover $\Theta_{\ep,\mu}$ verifies (\ref{pbthetal}), thence from 
(\ref{h2}) it follows that:
\begin{eqnarray}
     \eB{\lambda^{2}}|\Theta_{\ep,\mu}''(t)|^{2}&\leq &\frac{1}{\ep^{2}}
    |\Theta_{\ep,\mu}'(0)|^{2}e^{-2t/\ep} \eB{\lambda^{2}}\nonumber \\
    &\leq &  
    \frac{c_{3,\lambda}}{\ep^{2}} e^{-2t/\ep} 
    e^{2\lambda^{2}K_{4}\log(1+t)}.
    \label{thetamu2}
\end{eqnarray}
Using (\ref{umu2}), (\ref{thetamu2})  and (\ref{h2}) we get that
\begin{eqnarray*}
     \frac{ e^{2\lambda^{2}B_{\ep}(t)}}{b_{\ep}^{4}(t)}|U_{\ep,\mu}''(t) - 
       \Theta_{\ep,\mu}''(t)|^{2} &\leq   & \frac{2}{b_{\ep}^{4}(t)} 
       \eB{\lambda^{2}}(|U_{\ep,\mu}''(t)|^{2} + 
       |\Theta_{\ep,\mu}''(t)|^{2}) 
    \nonumber  \\
     & \leq & \frac{2c_{2,\lambda}}{\ep^{2} b_{\ep}^{2}(t)}(1+t)^{-2}  +
     \frac{2c_{3,\lambda}}{\ep^{2} b_{\ep}^{4}(t)}e^{-2t/\ep} e^{2\lambda^{2}K_{4}\log(1+t)}
    \nonumber  \\
     & \leq & \frac{1}{\ep^{2}}\left[c_{4,\lambda} + 
     c_{5,\lambda}(1+t)^{4}e^{-t}e^{2\lambda^{2}K_{4}\log(1+t)}\right] 
     %\nonumber \\
    \leq 
    \frac{ c_{6,\lambda}}{\ep^{2}}, 
    \label{diff1}
 \end{eqnarray*}
that is (\ref{S3}).
  \paragraph{\textsc{Estimate on 
  $V_{\ep,\lambda}''-\theta_{\ep,\lambda}''$}}
  Let $M=A_{[\lambda,\mu)}$ and $b(t) = b_{\ep}(t)$. Then 
  $V_{\ep,\lambda}$ and $\theta_{\ep,\lambda}$ are the solutions of the 
  corresponding problems (\ref{cpM}) and (\ref{cptheta}). Moreover since $A_{[\lambda,\mu)}$ is a bounded 
  operator we have that the related initial data $(v_{0},v_{1}) 
  \in D(M^{2})\times D(M^{3/2})$ and 
  $$|v_{0}|_{D(M^{2})}^{2} + |v_{1}|_{D(M^{3/2})}^{2} \leq 
  c_{7,\lambda}.$$ Let $\ep$ small in such a way that we can apply 
  Proposition \ref{p2} with these choices.
  Then since $w_{\ep}''= V_{\ep,\lambda}'' - 
       \theta_{\ep,\lambda}''$, and $\sigma_{M}^{2} = \lambda^{2}$, we 
       have that for every $t\geq 0$:
   \begin{equation}
       e^{2\lambda^{2}B_{\ep}(t)}\frac{1}{b_{\ep}^{4}(t)}|V_{\ep,\lambda}''(t) - 
       \theta_{\ep,\lambda}''(t)|^{2}\leq 
       c_{8,\lambda}.
       \label{S2}
   \end{equation}
 \paragraph{\textsc{Conclusion}}
 The inequality (\ref{D3}) in the general case is a straightforward consequence of 
 (\ref{S3}) and (\ref{S2}). \qed
 \subsection{A decomposition of $u_{\ep}$}
 Let $\uep$ be the solution of (\ref{pbm:h-eq}) - (\ref{pbm:h-data}) 
as in Theorem \ref{A} and let $u_{\ep,\nu}$ be defined as in (\ref{dec}).  Moreover let us set
 \begin{equation}
     |A^{1/2}\uep(t)|^{2}= \nu^{2}|u_{\ep,\nu}(t)|^{2} 
     +\alpha_{\ep,1}(t), \hspace{2em}  |A\uep(t)|^{2}= 
     \nu^{4}|u_{\ep,\nu}(t)|^{2} 
     +\alpha_{\ep,2}(t),  \label{Def}
 \end{equation}
  \begin{equation}
     \frac{ |\uep'(t)|^{2}}{b_{\ep}^{2}(t)}= 
    \frac{ |u_{\ep,\nu}'(t)|^{2} }{b_{\ep}^{2}(t)}
     +\alpha_{\ep,3}(t), \hspace{1em} \frac{ |A^{1/2}\uep'(t)|^{2}}{b_{\ep}^{2}(t)}= 
    \frac{ \nu^{2}|u_{\ep,\nu}'(t)|^{2} }{b_{\ep}^{2}(t)}
     +\alpha_{\ep,4}(t),
    \label{Defb}
 \end{equation}
 and
 \begin{equation}
     \eB{\nu^{2}} |u_{\ep,\nu}|^{2}= \beta_{\ep,0}(t), \hspace{1em} 
     \eB{\nu^{2}}\alpha_{\ep,1}(t) = \beta_{\ep,1}(t), \hspace{1em} 
     \eB{\nu^{2}}\alpha_{\ep,2}(t) = \beta_{\ep,2}(t),
     \label{Def1}
 \end{equation}
      \begin{equation}	  
      \eB{\nu^{2}}\alpha_{\ep,3}(t) = \beta_{\ep,3}(t),\hspace{2em} 
       \eB{\nu^{2}}\alpha_{\ep,4}(t) = \beta_{\ep,4}(t).
     \label{Def1b}
 \end{equation}
 In the proposition below we study the behaviour of 
 quantities defined in (\ref{Def1})-(\ref{Def1b}).
 \begin{prop} \label{p3}
   For $\ep$ small enough  the following properties hold true.
   \begin{enumerate}
       \item  For $t \rightarrow +\infty$ we have that:
      \begin{equation}
       \beta_{\ep,1}(t) \rightarrow 0, \hspace{2em}
         \beta_{\ep,2}(t) \rightarrow 0, \hspace{2em}
         \beta_{\ep,3}(t) \rightarrow 0,\hspace{2em}
	 \beta_{\ep,4}(t) \rightarrow 0,
       \label{lim2}
   \end{equation} 
   \begin{equation}
       \beta_{\ep,0}(t) \rightarrow L_{\ep} \in \re\setminus\{0\}.
       \label{lim1}
   \end{equation}
    \item  If $u_{0,\nu}\neq 0$ then there exists a constant 
       $K_{7}> 0$ independent of $\ep$ and $t$ such that
       \begin{equation} 
           \beta_{\ep,0}(t)\geq K_{7}, \hspace{2em} \forall t\geq 0.
           \label{stiun}
       \end{equation} 
       \end{enumerate}
   \end{prop}
 \paragraph{Proof of Proposition \ref{p3}}
  Let us denote by $c_{i}$ various constants that depend only on 
 $\nu$, $|u_{0}|_{D(A)}$ and $|u_{1}|_{D(A^{1/2})}$. 
  \subparagraph{\textsc{Proof of (\ref{lim2})}}
  Let us choose 
  $$\delta^{2}:= \nu^{2} + \frac{1}{K_{3}}.$$
 Let us assume that 
  $\ep$ is small enough so that we can use Theorem 
  \ref{t1} with  $\lambda = \delta$.

  We can rewrite the quantities in (\ref{Def}) and (\ref{Defb}) as:
  $$ \alpha_{\ep,h}(t) =  
 \sum_{k:\nu<\lambda_{k}<\delta}\lambda_{k}^{2h}|u_{\ep,k}(t)|^{2} + 
 |A^{h/2}U_{\ep,\delta}(t)|^{2} =  
 \alpha_{\ep,h,1}(t) +  \alpha_{\ep,h,2}(t), \hspace{1em} \mbox{for } 
 \,h = 1, 2,$$
  and for $h =3, 4$:
  $$\alpha_{\ep,h}(t) =  
 \frac{1}{b_{\ep}^{2}(t)}
 \left(\sum_{k:\nu<\lambda_{k}<\delta}\lambda_{k}^{2(h-3)}|u_{\ep,k}'(t)|^{2} + 
 |A^{(h-3)/2}U_{\ep,\delta}'(t)|^{2}\right) = 
 \alpha_{\ep,h,1}(t) +  \alpha_{\ep,h,2}(t).$$
 Since it holds true that:
$$ \frac{|A^{h/2}U_{\ep,\delta}'(t)|^{2}}{b_{\ep}^{2}(t)}
= \frac{\ep |A^{h/2}U_{\ep,\delta}'(t)|^{2}}{b_{\ep}(t)}
\frac{1}{\ep}\frac{1}{b_{\ep}(t)},
 $$ 
 thence thanks to (\ref{D1}) with $h=0$ and $h=1$ and 
 (\ref{h2}) we have that
 \begin{eqnarray*}
     &\displaystyle \eB{\nu^{2}}( \alpha_{\ep,1,2}(t) +  \alpha_{\ep,2,2}(t) + 
     \alpha_{\ep,3,2}(t)  + \alpha_{\ep,4,2}(t)) \leq  & \nonumber \\
    &\displaystyle   \leq\frac{ 
     c_{1}}{\ep}\frac{1}{b_{\ep}(t)} \eB{(\nu^{2}-\delta^{2})} 
     \leq  \frac{c_{2}}{\ep}(1+t) e^{-2\log(1+t)}.&
     \label{slim1}
 \end{eqnarray*}
 Hence we get that
 \begin{equation}
     \lim_{t\rightarrow + \infty} \eB{\nu^{2}}( \alpha_{\ep,1,2}(t) +  \alpha_{\ep,2,2}(t) + 
     \alpha_{\ep,3,2}(t) + \alpha_{\ep,4,2}(t) ) = 0.
     \label{slim2}
 \end{equation}
 
 For $\nu< \lambda_{k} < \delta$ let us now consider $b(t) = b_{\ep}(t)$ 
 and $M = A_{\{\lambda_{k}\}}$. 
 Then from (\ref{h2}) the function $b$ verifies (\ref{h2b}) and $M$ verifies (\ref{coercM}) with $\sigma_{M}^{2} = 
 \lambda_{k}^{2}$. Let $\ep$ small enough so that we can apply 
 Proposition \ref{p1} with such choices. We stress that since 
 $ \nu \leq \lambda_{k} \leq \delta$ we can take the smallness of $\ep$ independent 
 of $\lambda_{k}$. Since $\lambda_{k}> \nu$ and $\lambda_{k} 
 < \delta$ moreover from 
 (\ref{SL1b}) and (\ref{SL2b}) (with $h=1$) we have that:
 \begin{eqnarray}
     &\displaystyle\eB{\nu^{2}} ( \alpha_{\ep,1,1}(t) +  \alpha_{\ep,2,1}(t) + 
     \alpha_{\ep,3,1}(t) + \alpha_{\ep,4,1}(t) ) &\nonumber \\
     &\displaystyle \leq  c_{3} \sum_{k:\nu<\lambda_{k}<\delta}
     \left(2\lambda_{k}^{4} 
    |u_{\ep,k}(t)|^{2}+2 \frac{|u_{\ep,k}'(t)|^{2}}{b_{\ep}^{2}(t)}\right)\eB{\lambda_{k}^{2}}\eB{(\nu^{2}-\lambda_{k}^{2})}
    &\nonumber \\
   & \displaystyle \leq c_{4} 
     \displaystyle \sum_{k:\nu<\lambda_{k}<\delta}
     (\lambda_{k}^{2}+\lambda_{k}^{4} 
     +  1)(|u_{0,k}|^{2}+|u_{1,k}|^{2})\eB{(\nu^{2}-\lambda_{k}^{2})}& 
     \nonumber \\
     &\displaystyle \leq c_{5}  \displaystyle  \sum_{k:\nu<\lambda_{k}<\delta}(|u_{0,k}|^{2}+|u_{1,k}|^{2})
     \eB{(\nu^{2}-\lambda_{k}^{2})}.&
     \label{slim3}
 \end{eqnarray}
We can therefore passing to the limit in (\ref{slim3}) so that
\begin{eqnarray}
     &
    0 \leq \lim_{t\rightarrow +\infty}\eB{\nu^{2}} ( \alpha_{\ep,1,1}(t) +  \alpha_{\ep,2,1}(t) + 
     \alpha_{\ep,3,1}(t)+  \alpha_{\ep,4,1}(t)) \leq & \nonumber
     \\
     & \leq c_{5}  \displaystyle 
     \sum_{k:\nu<\lambda_{k}<\delta}\lim_{t\rightarrow +\infty}(|u_{0,k}|^{2}+|u_{1,k}|^{2})
     \eB{(\nu^{2}-\lambda_{k}^{2})} = 0. &
    \label{slim4}
\end{eqnarray}
From (\ref{slim2}) and ({\ref{slim4}) we get immediately (\ref{lim2}).
\subparagraph{\textsc{Proof of (\ref{lim1}) -  (\ref{stiun})}}
Let us set $y_{\ep}(t):= |u_{\ep,\nu}(t)|^{2}$. Then $y_{\ep}$ solves:
\begin{equation}
    y_{\ep}' = -2\nu^{2}b_{\ep}y_{\ep} -2\ep\langle u_{\ep,\nu}, 
    u_{\ep,\nu}''\rangle,
    \label{eqy}
\end{equation}
thence for all $t\geq 0$ we get that
\begin{equation}
    \eB{\nu^{2}}y_{\ep}(t) = |u_{0,\nu}|^{2} - 2\ep\int_{0}^{t}
    \langle u_{\ep,\nu}(s), 
    u_{\ep,\nu}''(s)\rangle e^{2\nu^{2}B_{\ep}(s)}ds.
    \label{eqyI}
\end{equation}
 Let us now estimate $\langle u_{\ep,\nu}, 
    u_{\ep,\nu}''\rangle$. Let us choose $M= A_{\{\nu\}}$, $b(t) =
    b_{\ep}(t)$ and let $\ep$ small so that we can apply Proposition 
    \ref{p1} and Proposition \ref{p2} (with $v_{\ep} = 
    u_{\ep,\nu}$).  This is possible since in such a case
    $$|u_{1,\nu}|_{D(M^{3/2})}^{2}+|u_{0,\nu}|_{D(M^{2})}^{2}  \leq 
    c_{6}( |u_{1,\nu}|^{2}+|u_{0,\nu}|^{2} ).$$
    Moreover clearly we have that
    $$|u_{\ep,\nu}(t)|^{2} = \nu^{-4}|Mu_{\ep,\nu}(t)|^{2}.$$
    Then from (\ref{SL1b})
    with $h=1$ and   (\ref{SL3b}), using  (\ref{h2}) (or 
    equivalently (\ref{h2b})) we obtain that
    \begin{eqnarray}
        |\langle u_{\ep,\nu}(t), 
    u_{\ep,\nu}''(t)\rangle|   \eB{\nu^{2}}& \leq & |u_{\ep,\nu}(t)|
    (|w_{\ep}''(t)| + |\theta_{\ep}''(t)|)\eB{\nu^{2}}
        \nonumber  \\
         & \leq & 
	 c_{7} ( \|u_{1,\nu}\|+|u_{0,\nu}| ) \left(\frac{|w_{\ep}''(t)|}{b_{\ep}^{2}(t)}b_{\ep}^{2}(t)
	 + \frac{1}{\ep}|\theta_{\ep}'(0)|
	 e^{-t/\ep} \right) e^{\nu^{2}B_{\ep}(t)} 
   \nonumber  \\
         & \leq &  c_{8} ( \|u_{1,\nu}\|^{2}+|u_{0,\nu}|^{2} ) 
	 \left(b_{\ep}^{2}(t) + \frac{1}{\ep}e^{-t/\ep} e^{\nu^{2}B_{\ep}(t)}\right)
        \nonumber  \\
         & \leq & c_{9}( \|u_{1,\nu}\|^{2}+|u_{0,\nu}|^{2} ) 
	 \left(\frac{1}{(1+t)^{2}} +  
	 \frac{1}{\ep} e^{-t/\ep}e^{\nu^{2}K_{4}\log(1+t)}\right).
	 \nonumber  \end{eqnarray}
	 Using that
	 $$\sup_{t\geq 0} e^{-t/2}e^{\nu^{2}K_{4}\log(1+t)} < + 
	 \infty,$$
	hence  we arrive at
        \begin{eqnarray}
        |\langle u_{\ep,\nu}(t), 
    u_{\ep,\nu}''(t)\rangle|   \eB{\nu^{2}} & \leq &
    c_{10}( \|u_{1,\nu}\|^{2}+|u_{0,\nu}|^{2} ) 
	 \left(\frac{1}{(1+t)^{2}} +  
	 \frac{1}{\ep} e^{-t/2\ep}\right).
        \label{stiprodotto}
    \end{eqnarray}
    
    From (\ref{stiprodotto}) thus we get for all $t\geq 0$ that
    \begin{eqnarray}
       \left| \int_{0}^{t}
    \langle u_{\ep,\nu}(s), 
    u_{\ep,\nu}''(s)\rangle e^{2\nu^{2}B_{\ep}(s)}ds\right| & \leq & \int_{0}^{t}
    |\langle u_{\ep,\nu}(s), 
    u_{\ep,\nu}''(s)\rangle| e^{2\nu^{2}B_{\ep}(s)}ds \nonumber \\
   & \leq  & c_{11}( \|u_{1,\nu}\|^{2}+|u_{0,\nu}|^{2} ),
        \label{stiint}
    \end{eqnarray}
    and also
   \begin{equation}\lim_{t\rightarrow +\infty}\int_{0}^{t}
    \langle u_{\ep,\nu}(s), 
    u_{\ep,\nu}''(s)\rangle e^{2\nu^{2}B_{\ep}(s)}ds = S_{\ep} \in \re.
     \label{limint}
    \end{equation}
 Therefore from (\ref{eqyI}) and (\ref{limint}) we have that there 
 exists
 $$\lim_{t\rightarrow + \infty} \eB{\nu^{2}}y_{\ep}(t) = 
 |u_{0,\nu}|^{2} - 2\ep S_{\ep}. $$
 We have  to prove that this limit is not zero.
 
 \subparagraph{Case $u_{0,\nu}\neq 0$.} 
By (\ref{stiint}), for $\ep$ small we 
have that
$$ 2\ep \left| \int_{0}^{t}
    \langle u_{\ep,\nu}(s), 
    u_{\ep,\nu}''(s)\rangle e^{2\nu^{2}B_{\ep}(s)}ds\right| \leq 
    \frac{1}{2} |u_{0,\nu}|^{2}, \quad\quad \forall t\geq 0,$$
    hence  (\ref{stiun}) follows from (\ref{eqyI}) and, as a 
    consequence, the limit in (\ref{lim1}) is different from zero.

    \subparagraph{Case $u_{0,\nu}= 0$.} Since $u_{1,\nu}\neq 0$, then 
    there exists a  single \emph{real} component of 
    $u_{1,\nu}$ different from zero, 
    that we indicate by $u_{1,\nu,r}$.   Let $u_{\ep,\nu,r}$ the 
    related component of $u_{\ep,\nu}$. 
    We prove that 
    \begin{equation}
        \lim_{t\rightarrow + \infty} \eB{\nu^{2}}|u_{\ep,\nu,r}(t)|^{2} 
	\neq 0.
        \label{lim5}
    \end{equation}
    This will be enough to prove that limit in (\ref{lim1}) is not 
    zero.
    
    To begin with, let us remark that there exists $T_{\ep}>0$ such that
    $$u_{\ep,\nu,r}'(T_{\ep}) = 0.$$
    Indeed if it is  not the case, then 
    $u_{\ep,\nu,r}$ is a strictly increasing or decreasing function 
    and since $u_{\ep,\nu,r}(0) = 0$, therefore we get that
    $$\lim_{t\rightarrow +\infty} u_{\ep,\nu,r} (t) \neq 0,$$
   but this is in contrast with (\ref{h1}), since $|A^{1/2}\uep(t)|^{2}\geq 
    \nu^{2}| u_{\ep,\nu,r} (t)|^{2}$ for all $t\geq 0$.
    
    Let now us set
    $$T_{\ep,0}:=\sup\{\tau\geq 0:\,  u_{\ep,\nu,r}' (t) \neq 0, \; 
    \forall t\in[0,\tau]\}.$$
  As seen before,  $T_{\ep,0}$ is a real positive number, 
  moreover 
  $$u_{\ep,\nu,r}' (T_{\ep,0}) = 0,$$
  and in $[0,T_{\ep,0}[$  the function
  $u_{\ep,\nu,r}$ is  strictly increasing or decreasing, so 
  that
  $$u_{\ep,\nu,r} (T_{\ep,0}) = P_{\ep}\neq  0.$$
  Therefore, as in (\ref{eqy}) - (\ref{eqyI}) for $t\geq T_{\ep,0}$ 
  we have that
  \begin{equation}
    e^{2\nu^{2}(B_{\ep}(t) - B_{\ep}(T_{\ep,0}))}|u_{\ep,\nu,r}(t)|^{2} = 
    P_{\ep}^{2} - 2\ep\int_{T_{\ep,0}}^{t}
     u_{\ep,\nu,r}(s) 
    u_{\ep,\nu,r}''(s)  e^{2\nu^{2}(B_{\ep}(s) - 
    B_{\ep}(T_{\ep,0}))}ds.
    \label{eqyI2}
\end{equation}
Now for $t\geq 0$, let us set
$ v_{\ep}(t) = u_{\ep,\nu,r}(t+ T_{\ep,0})$. Then $v_{\ep}$ verifies
(\ref{cpM}) with $M= A_{\{\nu\}}$ 
restricted to the single component $u_{\ep,\nu,r}$, $b(t) = 
b_{\ep}(t+T_{\ep,0})$  and initial data 
$v_{\ep}(0) = P_{\ep}$, $v_{\ep}'(0) = 0$. Thanks to (\ref{h2}) it is 
clear that the function $b$ verifies (\ref{h2b}). Therefore we can 
obtain as in (\ref{stiprodotto}) and (\ref{stiint}):
\begin{equation}\left|\int_{T_{\ep,0}}^{t}
     u_{\ep,\nu,r}(s) 
    u_{\ep,\nu,r}''(s)  e^{2\nu^{2}(B_{\ep}(s) - 
    B_{\ep}(T_{\ep,0}))}ds\right| \leq 
     c_{11}P_{\ep}^{2}, \quad \quad \forall t 
    \geq T_{\ep,0}.\label{stiint2}
\end{equation}
Only we have to specify that
$$\sup_{t\geq T_{\ep,0}} e^{-(t-T_{\ep,0})/2} 
e^{\nu^{2}K_{4}(\log(1+t) - \log(1+T_{\ep,0}))} \leq\sup_{t\geq 0} 
e^{-t/2} 
e^{\nu^{2}K_{4}\log(1+t)}  < +\infty.$$
Let now $\ep$ be small enough so that
$\ep c_{11} \leq 1/2,$
then from (\ref{eqyI2}) and (\ref{stiint2}) we get that
$$ e^{2\nu^{2}(B_{\ep}(t) - B_{\ep}(T_{\ep,0}))}|u_{\ep,\nu,r}(t)|^{2}  
\geq \frac{1}{2}P_{\ep}^{2},\quad\quad \forall t\geq T_{\ep,0},$$
thus the limit in (\ref{lim5}) is different from zero.
%% 
 % Since  for $t\geq T_{\ep,0}$ we have
 % $$
 %     e^{2\nu^{2}B_{\ep}(t)}|u_{\ep,\nu}(t)|^{2}\geq 
 %     e^{2\nu^{2}B_{\ep}(t)}|u_{\ep,\nu,r}(t)|^{2} 
 %     \geq \frac{1}{2}L_{\ep}^{2} e^{2\nu^{2}B_{\ep}(T_{\ep,0})}$$
 %     we have then shown that the limit in (\ref{lim1}) is different from 0.
 %%
 \qed
 \subsection{Proof of Theorem \ref{t2}}
 Let us assume that $\ep$ is small enough so that Theorem \ref{t1} 
 with $\lambda=\nu$  and 
 Proposition \ref{p3} hold true. 
  Let us moreover denote by $c_{i}$ various constants that depend only on 
 $\nu$, $|u_{0}|_{D(A)}$ and $|u_{1}|_{D(A^{1/2})}$ and by 
  $c_{i,\ep}$ constants that depend also on $\ep$.
  
  \paragraph{\textsc{Proof of (\ref{B1})}}
  Since the limit in (\ref{lim1}) is different from zero, hence there 
  exists $T_{\ep,1}\geq 0$ such that:
  $$\beta_{\ep,0}(t) \geq c_{1,\ep} > 0,\quad\quad \forall t\geq 
  T_{\ep,1}$$
  and in particular
  $$|u_{\ep,\nu}(t) |> 0,\hspace{2em} \forall t\geq T_{\ep,1}.$$
  Let us remark that if $u_{0,\nu}\neq 0$, then thanks to (\ref{stiun}) we 
  can take $T_{\ep,1} = 0$. 
  
 Thanks to (\ref{Def}) and (\ref{defb}) for $t\geq T_{\ep,1}$ we have 
 that
 \begin{eqnarray}
     b_{\ep}(t)\eB{\nu^{2}\gamma} & = & \left(\nu^{2}|u_{\ep,\nu}(t)|^{2} 
     +\alpha_{\ep,1}(t)\right)^{\gamma}\eB{\nu^{2}\gamma}
     \nonumber \\
      & = & \nu^{2\gamma}\beta_{\ep,0}^{\gamma}(t)
      \left(1 + 
      \frac{\alpha_{\ep,1}(t)}{\nu^{2}|u_{\ep,\nu}(t)|^{2}}\right)^{\gamma}
     \nonumber \\
      & = & \nu^{2\gamma}\beta_{\ep,0}^{\gamma}(t)\left(1 + 
      \frac{\beta_{\ep,1}(t)}{\nu^{2}\beta_{\ep,0}(t)}\right)^{\gamma}.
     \label{eqb1}
 \end{eqnarray}
 
Since for all $x\geq 0$ there exists $0\leq \xi\leq x$ such that
$$(1+x)^{\gamma}= 1 + \gamma (1+\xi)^{\gamma - 1}x,$$
then for $t\geq T_{\ep,1}$ we can rewrite (\ref{eqb1}) as
\begin{equation}
     b_{\ep}(t)\eB{\nu^{2}\gamma} = 
     \nu^{2\gamma}\beta_{\ep,0}^{\gamma}(t) \left(1 +\gamma (1+\xi)^{\gamma-1}
     \frac{\beta_{\ep,1}(t)}{\nu^{2}\beta_{\ep,0}(t)}\right)
     =:\nu^{2\gamma}\beta_{\ep,0}^{\gamma}(t) + \phi_{\ep}(t), 
    \label{eqb2}
\end{equation}
where if $\gamma < 1$ then
\begin{equation}
    0\leq  \phi_{\ep}(t)\leq c_{1} \beta_{\ep,0}^{\gamma-1}(t)
  \beta_{\ep,1}(t),
    \label{stif1}
\end{equation}
while if $\gamma \geq 1$ then 
\begin{equation}
    0\leq  \phi_{\ep}(t)\leq c_{2} \beta_{\ep,0}^{\gamma}(t)
    \left(1 + 
      \frac{\beta_{\ep,1}(t)}{\beta_{\ep,0}(t)}\right)^{\gamma-1}
      \frac{\beta_{\ep,1}(t)}{\beta_{\ep,0}(t)}\leq  
      c_{2} \left(\beta_{\ep,0}(t) +
      \beta_{\ep,1}(t)\right)^{\gamma-1}
      \beta_{\ep,1}(t).
    \label{stif2}
\end{equation}

Integrating (\ref{eqb2}) we get that
\begin{equation}
\eB{\nu^{2}\gamma}  - e^{2\nu^{2}\gamma B_{\ep}(T_{\ep,1})} = 
2\nu^{2}\gamma \left[\int_{T_{\ep,1}}^{t} 
\nu^{2\gamma}\beta_{\ep,0}^{\gamma}(s) + \phi_{\ep}(s) \, ds\right],\quad\quad 
\forall t\geq T_{\ep,1}.  
    \label{eqB}
\end{equation}
From (\ref{lim1}), (\ref{lim2}) and (\ref{stif1}) or (\ref{stif2}) we 
immediately obtain that
\begin{eqnarray}
    \lim_{t\rightarrow + \infty}\displaystyle \frac{1}{1+t}\int_{T_{1}}^{t} 
\beta_{\ep,0}^{\gamma}(s) \, ds & = & L_{\ep}^{\gamma},
    \label{L1}  \\
    \lim_{t\rightarrow + \infty} \phi_{\ep}(t) & = & 0 ,\label{Lfi}  \\
     \lim_{t\rightarrow + \infty}\displaystyle \frac{1}{1+t}\int_{T_{1}}^{t} 
\phi_{\ep}(s) \, ds & = & 0.
    \label{L2}  
\end{eqnarray}
from (\ref{eqB}), (\ref{L1}) and (\ref{L2}) we then deduce that
\begin{equation}
   \lim_{t\rightarrow + \infty} \frac{1}{1+t}\eB{\nu^{2}\gamma} = 
    2\nu^{2(\gamma+1)}\gamma L_{\ep}^{\gamma}.
   \label{L3}
\end{equation}
This non zero limit proves (\ref{B1}) with constants depending on $\ep$.

Let us now assume that $u_{0,\nu}\neq 0$ so that $T_{\ep,1}= 0$ and 
(\ref{stiun}) holds true. Then from (\ref{eqB})
we obtain that
\begin{equation}
   \eB{\nu^{2}\gamma} \geq 1 +  2\nu^{2(\gamma+1)}\gamma 
   K_{7}^{\gamma}t, \hspace{2em} 
   \forall t\geq 0.
    \label{SB1}
\end{equation}
Moreover since $\uep = U_{\ep,\nu}$, from (\ref{D1}) of Theorem 
\ref{t1} (with $h = 0$ and $\lambda = \nu$) we get that
%\begin{equation}
   $$ b_{\ep}(t)\eB{\nu^{2}\gamma} = (\eB{\nu^{2}}|A^{1/2}u_{\ep}(t)|^{2})^{\gamma} \leq 
    \gamma_{0,\nu}^{\gamma} $$
%\end{equation}
hence for all $t\geq 0$ we have that
\begin{equation}
    \eB{\nu^{2}\gamma} = 1 +2\nu^{2}\gamma \int_{0}^{t} 
    (e^{2\nu^{2}B_{\ep}(s)}|A^{1/2}u_{\ep}(s)|^{2})^{\gamma}ds \leq 
    1 +2\nu^{2}\gamma\,
     \gamma_{0,\nu}^{\gamma} t.
    \label{SB2}
\end{equation}
Thus from (\ref{SB1}) and (\ref{SB2}) we obtain (\ref{B1}) with constants 
independent of $\ep$.
 \paragraph{\textsc{Proof of (\ref{B2})}}
From (\ref{eqB}), for $t\geq T_{\ep,1}$
we have that
\begin{equation}
    B_{\ep}(t) = \frac{1}{2\nu^{2}\gamma}\log\left(e^{2\nu^{2}\gamma 
    B_{\ep}(T_{\ep,1})} +2\nu^{2}\gamma \int_{T_{\ep,1}}^{t} 
\nu^{2\gamma}\beta_{\ep,0}^{\gamma}(s) + \phi_{\ep}(s) \, ds\right) .
     \label{eqb3}
\end{equation}
Taking the derivative of (\ref{eqb3}) we obtain that
%\begin{equation}
  $$  b_{\ep}(t) = \frac{\nu^{2\gamma}\beta_{\ep,0}^{\gamma}(t) + \phi_{\ep}(t)}{e^{2\nu^{2}\gamma 
    B_{\ep}(T_{\ep,1})} +\displaystyle 2\nu^{2}\gamma \int_{T_{\ep,1}}^{t} 
\nu^{2\gamma}\beta_{\ep,0}^{\gamma}(s) + \phi_{\ep}(s) \, ds}.
   $$
   %\label{eqb4}
%\end{equation}
Using (\ref{L1}), (\ref{Lfi}), (\ref{L2}) and (\ref{lim1}) we get that
$$\lim_{t\rightarrow +\infty} (1+t)b_{\ep}(t) = \lim_{t\rightarrow +\infty}
\frac{\nu^{2\gamma}\beta_{\ep,0}^{\gamma}(t) + \phi_{\ep}(t)}{\displaystyle\frac{1}{1+t}
\left[e^{2\nu^{2}\gamma 
    B_{\ep}(T_{\ep,1})} + 2\nu^{2}\gamma \int_{T_{\ep,1}}^{t} 
\nu^{2\gamma}\beta_{\ep,0}^{\gamma}(s) + \phi_{\ep}(s) \, ds\right]} = 
\frac{1}{2\nu^{2}\gamma},$$
that is (\ref{B2}).
\paragraph{\textsc{Limit of $|A^{1/2}\uep|$}}
We prove that
\begin{equation}
        \lim_{t\rightarrow+\infty}(1+t)^{1/\gamma}|A^{1/2}\uep(t)|^{2} = 
	\frac{1}{(2\nu^{2}\gamma)^{1/\gamma}}.
        \label{B3}
    \end{equation}
To this end it is enough to remark that 
$$(1+t)^{1/\gamma}|A^{1/2}u_{\ep}(t)|^{2} = 
((1+t)b_{\ep}(t))^{1/\gamma}$$
and use (\ref{B2}).
\paragraph{\textsc{Proof of (\ref{B31b})}}
From (\ref{Def}) we have that 
\begin{eqnarray}
    (1+t)^{1/\gamma}\nu^{2}|u_{\ep,\nu}(t)|^{2} & = &  (1+t)^{1/\gamma}
    |A^{1/2}u_{\ep}(t)|^{2} - (1+t)^{1/\gamma}\alpha_{\ep,1}(t)
    \nonumber \\
     & = & (1+t)^{1/\gamma}
    |A^{1/2}u_{\ep}(t)|^{2} - 
    (1+t)^{1/\gamma}e^{-2\nu^{2}B_{\ep}(t)}\beta_{\ep,1}(t).
    \label{eqb5}
\end{eqnarray}
From (\ref{B1}) we know that 
\begin{equation}
    (1+t)^{1/\gamma}e^{-2\nu^{2}B_{\ep}(t)}\leq c_{2,\ep}, \quad 
    \quad \forall t\geq 0,
    \label{limitazione}
\end{equation}
 hence from 
(\ref{B3}),  (\ref{lim2}) and (\ref{eqb5}) we deduce that
\begin{equation}
  \lim_{t\rightarrow 
+\infty}(1+t)^{1/\gamma}\nu^{2}|u_{\ep,\nu}(t)|^{2}  = 
\frac{1}{(2\nu^{2}\gamma)^{1/\gamma}},
    \label{limun}
\end{equation}
whence (\ref{B31b}) immediately follows.

\paragraph{\textsc{Proof of (\ref{B5})}} 

Thanks to (\ref{D3}) in Theorem \ref{t1} we have for all $t\geq 0$ 
that:
\begin{equation}
    |U_{\ep,\nu}''(t)| \leq 
    |U_{\ep,\nu}''(t) - 
    \Theta_{\ep,\nu}''(t)| + | \Theta_{\ep,\nu}''(t)|\leq
   \sqrt{\gamma_{\ep,\nu}}\,  b_{\ep}^{2}(t)e^{-\nu^{2}B_{\ep}(t)} + 
     \frac{1}{\ep}|\Theta_{\ep,\nu}'(0)|e^{-t/\ep}.
    \label{stideru0} 
\end{equation}
%% 
 % hence 
 % \begin{equation}
 %    |U_{\ep,\nu}''(t)|  \leq 
 %      \sqrt{\gamma_{\ep,\nu}}\,  b_{\ep}^{2}(t)e^{-\nu^{2}B_{\ep}(t)} + 
 %      \frac{1}{\ep}|\Theta_{\ep,\nu}(0)|e^{-t/\ep}, \quad\quad \forall 
 %      t\geq 0.
 %     \label{stideru}
 % \end{equation}
 %%
Using (\ref{limitazione}) and (\ref{h2}) in (\ref{stideru0}) hence we get that
\begin{equation}
     |U_{\ep,\nu}''(t)|^{2} \leq 
    c_{3,\ep}\left(\frac{1}{(1+t)^{2+1/(2\gamma)}}+ e^{-t}\right)^{2} 
    \leq c_{4,\ep} \frac{1}{(1+t)^{4+1/\gamma}},
    \label{stideru2}
\end{equation}
that is (\ref{B5}), since $\uep= U_{\ep,\nu}$. 

 \paragraph{\textsc{Existence of $u_{\ep,\infty}$}}
Thanks to  (\ref{coerc}), (\ref{limitazione}), (\ref{Def}), 
(\ref{Def1}) and (\ref{h2}) we have that for 
all $t\geq 0$:
$$(1+t)^{1/\gamma} |\uep(t) - u_{\ep,\nu}(t)|_{D(A)}^{2}
\leq c_{3}(1+t)^{1/\gamma} |A \uep(t) - A u_{\ep,\nu}(t)|^{2} \leq
c_{4,\ep} \beta_{\ep,2}(t),$$
$$(1+t)^{2+1/\gamma} |\uep'(t) - u_{\ep,\nu}'(t)|_{D(A^{1/2})}^{2} 
\leq c_{5,\ep}(1+t)^{2}b_{\ep}^{2}(t) \beta_{\ep,4}(t) \leq 
c_{6,\ep}\beta_{\ep,4}(t) ,$$
hence from (\ref{lim2}) we obtain that
$$\lim_{t\rightarrow + \infty}(1+t)^{1/\gamma} 
|\uep(t) - u_{\ep,\nu}(t)|_{D(A)}^{2} +
(1+t)^{2+1/\gamma} |\uep'(t) - u_{\ep,\nu}'(t)|_{D(A^{1/2})}^{2} = 0.$$
Therefore for proving  (\ref{LIM}) we have only to show that the 
functions 
$(1+t)^{1/(2\gamma)}u_{\ep,\nu}(t)$ and 
$(1+t)^{1+1/(2\gamma)}u_{\ep,\nu}'(t)$ have the required limits. 
Since 
$$u_{\ep,\nu}'(t) = -\nu^{2}b_{\ep}(t) u_{\ep,\nu}(t) - \ep 
u_{\ep,\nu}''(t),$$
then  we have that
$$e^{\nu^{2}B_{\ep}(t)}u_{\ep,\nu}(t) =
u_{0,\nu} - \ep \int_{0}^{t} 
e^{\nu^{2}B_{\ep}(s)}u_{\ep,\nu}''(s)\, ds.$$
Thanks to (\ref{B5}) it is clear that for all $t\geq 0$:
\begin{equation}
    |u_{\ep,\nu}''(t)| \leq |\uep''(t)| \leq 
   \sqrt{ K_{\ep}}\frac{1}{(1+t)^{2+1/(2\gamma)}},
    \label{stiuder}
\end{equation}
thus using once again (\ref{B1}) we obtain that there exists 
$$\lim_{t\rightarrow 
+\infty}\int_{0}^{t}e^{\nu^{2}B_{\ep}(s)}u_{\ep,\nu}''(s)\, ds = 
\alpha_{\ep,\nu}\in H_{\{\nu\}}.$$ 
Applying (\ref{L3}) we  finally arrive at
\begin{eqnarray}
    \lim_{t\rightarrow+\infty}(1+t)^{1/(2\gamma)} u_{\ep,\nu}(t)& = &   
    \lim_{t\rightarrow+\infty}(1+t)^{1/(2\gamma)}
    e^{-\nu^{2}B_{\ep}(t)}e^{\nu^{2}B_{\ep}(t)} 
    u_{\ep,\nu}(t)\nonumber \\
     & = & \frac{1}{\left(2\nu^{2(\gamma+1)}\gamma 
     L_{\ep}^{\gamma}\right)^{1/(2\gamma)}}(u_{0,\nu} - \ep
     \alpha_{\ep,\nu}) = u_{\ep,\infty} \in H_{\{\nu\}} .\label{limuA}
\end{eqnarray}
Furthermore we have also that
$$(1+t)^{1+1/(2\gamma)}u_{\ep,\nu}'(t) = 
-\nu^{2}(1+t)b_{\ep}(t)(1+t)^{1/(2\gamma)}u_{\ep,\nu}(t) - \ep 
(1+t)^{1+1/(2\gamma)}u_{\ep,\nu}''(t),$$
therefore from (\ref{B2}), (\ref{limuA}) and (\ref{stiuder}) we get 
that
\begin{equation}
    \lim_{t\rightarrow+\infty}(1+t)^{1+1/(2\gamma)}u_{\ep,\nu}'(t)  = 
    -\frac{1}{2\gamma}u_{\ep,\infty} \in H_{\{\nu\}}.
    \label{limderuA}
\end{equation}
From (\ref{limuA}), (\ref{limderuA}) and (\ref{B31b}) 
the existence of the required non zero limits immediately follows.
\qed

\label{NumeroPagine}

\end{document}